\documentclass[11pt]{article}
\usepackage{fullpage,amsmath,amssymb,graphicx,tikz,caption,color,enumitem}
\usetikzlibrary{backgrounds, positioning}
\numberwithin{figure}{section}
\newtheorem{theorem}{Theorem}[section]
\newtheorem{lemma}[theorem]{Lemma}
\newtheorem{corollary}[theorem]{Corollary}
\newtheorem{observation}[theorem]{Observation}
\newcommand{\del}{\backslash}
\newcommand{\ov}{\overline}
\usepackage{listings,xcolor}
\begin{document}

\begin{center}
{\Large\bf On extending induced subgraphs}

\medskip
Guoli Ding\footnotemark[1],
Jake Weber\footnotemark[2],
and Weihao Xia\footnotemark[3]\\
Department of Mathematics, Louisiana State University,
Baton Rouge, LA 70803, USA
\end{center}

\footnotetext[1]{Email: ding@math.lsu.edu}
\footnotetext[2]{Email: jwebe33@lsu.edu}
\footnotetext[3]{Email: wxia3@lsu.edu}

\begin{abstract}
When characterizing graphs defined by forbidden induced subgraphs, it is often desirable to understand how a graph $G$ can be generated incrementally from a given induced subgraph $H$. In this paper, we study a number of problems within this framework. First, for several notions of connectivity, we determine all minimal induced subgraphs of $G$ that contain $H$ and attain the connectivity of $G$. Second, we establish how to grow $H$ inside $G$ while preserving twin-freeness and certain forms of connectivity. Third, we establish {\it finite test theorems} for various graph classes, proving that agreement with the recursively generated graphs for a few initial iterations forces agreement at all subsequent iterations. Consequently, such an infinite graph family can be completely determined by finitely many of its initial generations. As an application, we characterize $\{\mathbb L_1,\mathbb L_3\}$-free graphs, where $\mathbb L_1$ is the claw and $\mathbb L_3$ is another forbidden graph for line graphs. Our structural results reduce this problem to a finite verification, which is completed through computer-assisted computation.
\end{abstract}

\section{Introduction}

Induced subgraphs play a fundamental role in structural graph theory. Characterizing   graph properties through forbidden induced subgraphs has led to major breakthroughs, most notably the strong perfect graph theorem. Conversely, graph classes defined by excluded induced subgraphs often exhibit deep and elegant structures, as demonstrated by the structural characterization of claw-free graphs \cite{claw-free}. Research in both directions remains highly active today, driven in part by major open problems such as the Gy\'arf\'as-Sumner conjecture and the Erd\"os-Hajnal conjecture.

To study a graph class $\cal G$ closed under a given notion of graph containment, one often needs to analyze how a general graph $G\in\cal G$ is constructed from a contained graph $H$ via local extensions within $\cal G$. This perspective has been remarkably successful in graph minor theory, where deep structural results are established using various splitter theorems. In the context of induced subgraphs, analogous ideas provide a natural framework for structural characterization and recursive generation. In this paper, we investigate several problems arising from this perspective.

We use standard terminology and notation from graph theory. For definitions not explicitly provided here, we generally follow Diestel \cite{diestel}. For graphs $H$ and $G$, we write $H\subseteq G$ if $H$ is an induced subgraph of $G$, we write $H\preceq G$ if $H$ is isomorphic to an induced subgraph of $G$, and we write $H\subset G$ or $H\prec G$, respectively, if additionally $|H|<|G|$. For a set $\cal F$ of graphs, we say that $G$ is $\cal F$-{\it free} if no graph $F\in\cal F$ satisfies $F\preceq G$. 

As connectivity plays a central role throughout the paper, we formalize several connectivity concepts relevant to our results. A graph $G$ is {\it trivial} if $|G|=1$, and its {\it complement} is denoted by $\ov G$. A nontrivial graph $G$ is {\it doubly connected} if both $G$ and $\ov G$ are connected. A vertex $v$ of $G$ is {\it simplicial} if $N_G(v)$ is a clique. Let $G_1,G_2$ be disjoint connected nontrivial graphs, and let $v_i$ be a simplicial vertex of $G_i$ for $i=1,2$. Then 
\begin{enumerate}
\parskip -4pt
\item[$\bullet$] a {\it simplicial sum} of $G_1,G_2$ (over $v_1,v_2$) is obtained by identifying $v_1$ with $v_2$;
\item[$\bullet$] a {\it simplicial join} of $G_1$, $G_2$ (over $v_1$, $v_2$), where $d_{G_i}(v_i)\ge2$ for $i=1,2$, is obtained from $(G_1- v_1)\cup (G_2- v_2)$ by adding all edges between $N_{G_1}(v_1)$ and $N_{G_2}(v_2)$.
\end{enumerate}
If $G$ is a simplicial-sum or -join of $G_1,G_2$, it is easy to see that $G_1, G_2\prec G$. A graph is {\it reducible} if it is disconnected or can be expressed as a simplicial-sum or -join. A nontrivial graph is {\it irreducible} if it is not reducible. Reducibility can be viewed as a form of low connectivity: every reducible graph can be reconstructed from two smaller graphs. We also note that being reducible is essentially equivalent to admitting a 1-{\it join}, a concept introduced in \cite{claw-free}. In general, every claw-free graph admitting a 1-join is reducible, but the converse is not true, as illustrated by the simplicial sum of two copies of $K_5$.

Let $X,Y\subseteq V(G)$ be non-empty and disjoint. We say that $X$ is {\it complete}, {\it anticomplete}, or {\it mixed} to $Y$ if all, none, or some (but not all) possible edges between $X,Y$ are present in $G$, respectively. When $X$ or $Y$ is a singleton, we may identify it with its unique element. We also apply this terminology to subgraphs by identifying each subgraph with its vertex set.

When constructing a graph $G$ from an induced subgraph $H$, the initial step frequently involves augmenting $H$ inside $G$ to match the connectivity of $G$. However, the resulting induced subgraph can easily become excessively large, limiting the power of this approach. Indeed, a minimal induced subgraph of $G$ that contains $H$ and satisfies the desired connectivity condition can have arbitrarily large order relative to $|H|$. To illustrate this, let $H$ be the disjoint union of two triangles, and let $G$ be obtained by connecting the triangles with a path. Then $G$ is connected and is itself the minimal connected induced subgraph of $G$ containing $H$, yet $|G| - |H|$ can be made arbitrarily large independently of $|H|$. 

In Section 2, we resolve this issue in several key special cases. For instance, when both $G$ and $H$ are doubly connected and $G$ contains a vertex complete to $H$, we completely determine all minimal doubly connected induced subgraphs of $G$ containing $H$ and a vertex complete to $H$. We also resolve the problem when $G$ is irreducible while $H$ is a simplicial sum.

Besides connectivity considerations, structural simplification can also be achieved by eliminating equivalent vertices. Let $N_G[v]$ (or $N[v]$) denote the {\it closed neighborhood} $N_G(v) \cup \{v\}$ of a vertex $v$. Distinct vertices $u$ and $v$ are {\it twins} if $N_G[u] = N_G[v]$. Since twins play identical structural roles, a graph can be simplified by deleting one vertex of the pair. A graph is {\it prime} if it contains no twins. We note that while ``prime" has other definitions in the literature, we use it here exclusively to mean twin-free. 

In Section 3, we study the recursive generation of prime graphs under different connectivity settings. We establish several splitter theorems, including the following two central results. For any two connected prime graphs $H$ and $G$ satisfying $H\prec G$, there exists a sequence $H\cong G_0 \prec G_1\prec \cdots \prec G_k\cong G$ of connected prime graphs with $|G_i|-|G_{i-1}|=1$ for $i=1,\ldots, k$. Similarly, if $H$ and $G$ are doubly connected prime graphs satisfying $H\prec G$, and $G$ is not on a specific list of exceptional graphs, then there exists a sequence $H\cong G_0 \prec G_1\prec \cdots \prec G_k\cong G$ of doubly connected prime graphs with $|G_i|-|G_{i-1}|\le 4$ for $i=1, ..., k$. These results mirror the role of classical splitter theorems in graph minor theory and provide a robust framework for structural induction.

In Section 4, we establish three {\it finite test theorems}. The splitter theorems from Section 3 enable the iterative generation of prime graphs while preserving prescribed structural properties. Suppose that, for some sufficiently large integer $k$, the generated graphs of order at most $k$ perfectly coincide with the members of a graph class $\cal G$ up to that same order. A fundamental question is whether this finite observation guarantees that all subsequently generated graphs must also coincide exactly with $\cal G$. Although this implication fails in general, we identify two families of graph classes for which the implication does hold. These results provide a definitive termination criterion for the recursive generation process. As far as we know, these are the first results of this type in structural graph theory.

Finally, in Section 5, we characterize $\{\mathbb L_1,\mathbb L_3\}$-free graphs, where $\mathbb L_1$ (the claw) and $\mathbb L_3$ are two forbidden graphs for line graphs (see Figure~\ref{fig:line7}). Our proof is self-contained and does not rely on the general structure theorem for claw-free graphs \cite{claw-free}. Instead, we first apply standard reductions to reduce the problem to the characterization of $\{\mathbb L_1,\mathbb L_3\}$-free irreducible prime graphs. The structural results developed in the previous sections then reduce the remaining task to a finite computer-assisted verification. This application illustrates how recursive generation and finite test theorems can reduce a structural characterization problem to finite computation, thereby cleanly separating the theoretical and computational components of the proof.

We close this section with a few more definitions and observations. For any positive integer $k$, let $[k]=\{1,...,k\}$. If $H$ is a subgraph of $G$, then we denote $G-V(H)$ by $G-H$ and $G[V(H)]$ by $G[H]$. Similarly, if $H\subseteq G$ and $X\subseteq V(G-H)$, then we denote $G[V(H)\cup X]$ by $H+X$. 
The following observations are used explicitly and implicitly throughout the paper. We omit the proofs since they are straightforward.

\begin{observation}\label{obs:basic}
Let $G=(V,E)$ be a graph. Then the following hold.
\begin{enumerate}[label=(\arabic*)]
\parskip -4pt
\item If $x,y$ are twins in $G$, then they are also twins in $G-z$ for any $z\in V-x-y$. 
\item If $G$ is not prime but $G-x$ is prime for some $x\in V$, then $x$ has a twin in $G$.
\item If $G$ is connected, claw-free, and has a cut vertex, then $G$ is a simplicial sum.
\item Let $G$ be the simplicial join of $G_1$ and $G_2$ over $v_1$ and $v_2$. If $H\subseteq G$, where $H$ is claw-free and $H\cap (G_1-N_{G_1}[v_1])\ne \emptyset$ $\ne H\cap (G_2-N_{G_2}[v_2])$, then $H$ is reducible.
\item If $G$ is claw-free and prime, and $H$ is a component of $\ov G$ satisfying $3\le |H|<|G|$, then $H$ is doubly connected.
\end{enumerate}
\end{observation}

\section{Connectivity remediation}

In this section, we consider, under three different settings of connectivity, how to enhance the connectivity of an induced subgraph $H$ within a better connected host graph $G$. More specifically, we are interested in determining all minimal graphs $J$ satisfying $H\preceq J \preceq G$ such that $J$ inherits the connectivity of $G$. We solve this problem in three key special cases.

When $G$ is connected while $H$ consists of two components, constructing a desired graph $J$ is straightforward. One may simply add a shortest path in $G$ joining the two components of $H$. However, as illustrated in the introduction, the resulting graph $J$ can be arbitrarily larger than $H$. Such unbounded growth is inconsistent with our goal of improving $H$ through local operations. Nevertheless, the following simple observation shows that this issue can be avoided when at least one of the two components of $H$ is a path. 

\begin{theorem}\label{thm:pathext}
Let $G$ be connected and let $H\prec G$. If $H$ is disconnected and has a component $P$ that is a nontrivial path, then there exists $J$ such that $H\prec J\subseteq G$, $|J| - |H| \le |P|-1$, and $J$ has fewer components than $H$.
\end{theorem}

\noindent{\bf Proof.} Without loss of generality, we assume $H\subset G$. Let us also assume that $H$ is the disjoint union of $H_1$ and $P$. Since $G$ is connected, it has a shortest path $Q=v_1v_2\cdots v_k$ between $H_1$ and $P$, where $v_1\in V(H_1)$. If $|Q-H|=k-2\le |P|-1$, then $J=G[H\cup Q]$ satisfies the requirements. On the other hand, if $k-2>|P|-1$, it allows us to define $Q'=v_3Qv_t\cong P$ and $Q''=v_1Qv_t$, where $t=|P|+2$. Then $G[H_1\cup Q']\cong H$ and $|H_1\cup Q''|=|H|+1\le |H|+|P|-1$. These clearly imply that $J=G[H_1\cup Q'']$ satisfies the requirements. \hfill$\Box$\medskip

It is worth noting that, if $P$ is a trivial path, then the graph $J$ constructed in the proof still satisfies all the conclusions of the theorem, except that the bound is replaced by $|J|-|H|\le |P|$.

Next, we consider the case when $G$ is doubly connected. It turns out that this stronger form of connectivity is much harder to achieve than the ordinary connectivity. To explain the main obstacle, we assume, without loss of generality, that $H$ is disconnected; otherwise, $\ov H$ is disconnected, and we may work in $\ov G$ instead. Even when $H$ has only two components and $J$ is allowed to be arbitrarily large, the strategy of joining the components of $H$ by a shortest path may still fail. Indeed, suppose that such a path has a unique internal vertex $u$ and that $u$ is complete to $H$. Then adding $u$ to $H$ makes the graph connected, yielding a graph $H_1$, but at the same time $\ov{H_1}$ is disconnected. One might then try to restore connectivity in the complement by joining its components in $\ov G$, but the same phenomenon may recur: repairing connectivity on one side may simultaneously destroy it on the other.

Breaking this cycle is the key difficulty in extending $H$ to a doubly connected graph. The following problem captures the essential obstruction and leads to our second result of this section. Suppose $H$ is doubly connected and $G-H$ has a vertex that is complete to $H$. What are the minimal graphs that contain $H$ and a vertex complete to $H$, while also being doubly connected?  

\begin{figure}[ht]
\begin{center}
\includegraphics[scale=0.7]{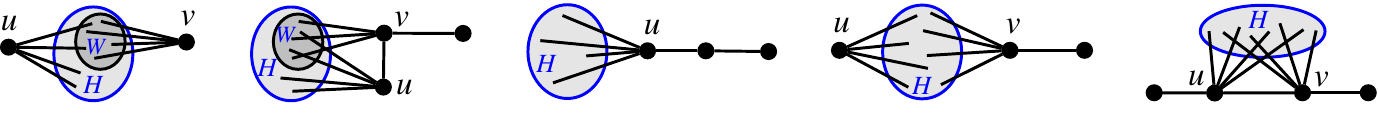}
\caption{$H$-remedies and $P_4$-augmentations of $H$.}
\label{fig:remedy}
\end{center}
\end{figure}

To answer this question, we introduce two types of constructions. A $P_4$-{\it augmentation} of $H$ is any of the last three graphs shown in Figure \ref{fig:remedy}, where $u$ and $v$ are complete to $H$. An $H$-{\it remedy} is either of the first two graphs in Figure \ref{fig:remedy}, where $\emptyset \ne W \subset V(H)$, $u$ is complete to $H$, and the neighborhood of $v$ in $H$ is exactly $W$. It is straightforward to verify that the last four graphs in Figure \ref{fig:remedy} are doubly connected, as is the first, provided that $\ov H$ is connected. Our next theorem characterizes precisely the minimal graphs that answer the above question.

\begin{theorem}\label{thm:remedy}
Let $G$ be doubly connected and let $H\subset G$. If some vertex of $G-H$ is complete to $H$, then there exists $J\subseteq G$ such that $J$ is either an $H$-remedy or a $P_4$-augmentation of $H$.
\end{theorem}

\noindent{\bf Proof.} Let $V=V(G)$, $X=V(H)$, and $k$ be the largest nonnegative integer satisfying the following. There exist disjoint subsets $Y = \{y_1, ... , y_s\}$ and $Z=\{z_1,...,z_t\}$ of $V\backslash X$, where $s =\lceil k/2\rceil$ and $t=\lfloor k/2\rfloor$ such that the following conditions are satisfied in some $\widetilde G\in\{G,\ov G\}$:\\  
\indent a. $Y$ is a clique, $Z$ is a stable set, and $y_iz_j$ is an edge if and only if $i\le j$;\\ 
\indent b. $X$ is complete to $Y$ and anticomplete to $Z$;\\ 
\indent c. $V\backslash U$, where $U=X\cup Y\cup Z$, has a vertex $w$ anticomplete to $U$.\\ 
Our assumption on $H$ implies that (a-c) hold trivially with  $\widetilde G=\ov G$, $k=0$, and $Y=Z=\emptyset$. This observation ensures the existence of $k$. In general, for the largest $k$, the choice of $(\widetilde G, Y, Z)$ may not be unique. Our proof does not depend on which $(\widetilde G, Y, Z)$ is chosen, except for the case $k=0$, where we require $(\widetilde G, Y, Z)$ to be $(\ov G,\emptyset, \emptyset)$. Since $\widetilde G$ is connected, it has a shortest path between $w$ and $U$. By replacing $w$ with a vertex on the path, if necessary, we assume that the path is $uvw$, where $u\in U$. We distinguish among three cases depending on how $v$ is adjacent to $X$ in $\widetilde G$.

We first assume that $v$ is mixed to $X$. If $\widetilde G =\ov G$, then $G[X\cup \{v,w\}]$ is the first graph in the figure. So we assume $\widetilde G = G$. According to how we choose $(\widetilde G, Y, Z)$ we conclude that $k>0$ and thus $Y\ne\emptyset$. If $vy_1\not\in E(G)$ then $G[X\cup \{v,y_1\}]$ is the first graph; if $vy_1\in E(G)$ then $G[X\cup \{v,w,y_1\}]$ is the second graph.

Next, assume that $v$ is anticomplete to $X$. If $v$ is adjacent to some $y_i$, then $\widetilde G[X\cup \{y_i,v,w\}]$ is the third graph in the figure. If $v$ is anticomplete to $Y$, then $v$ is adjacent to some $z_i$, and thus $\widetilde G[X\cup \{z_i,y_i,v\}]$ is the third graph. In both cases, $G$ contains the third or the fourth graph.

Finally, $v$ is complete to $X$. We claim that $v$ is not complete to $U$.
Otherwise, let $G'$ be the complement of $\widetilde G$, and set $Y'= Z\cup\{w\}$, $Z'=Y$, and $w'=v$. Then, in $G'$, $Y'$ is a clique, $Z'$ is a stable set, and conditions (b) and (c) hold. Moreover, if $w$ is indexed first in $Y'$, the adjacency pattern required in (a) holds for $k+1$. Hence $G',Y',Z'$, and $w'$ satisfy (a)-(c) for $k+1$, contradicting the maximality of $k$, which proves the claim. This claim implies that $v$ is not complete to $Y\cup Z$. If $v$ is nonadjacent to some $y_i$, then $\widetilde G[X\cup\{y_i,v,w\}]$ is the fourth graph in the figure, which implies that $G$ contains the third or the fourth graph. So, $v$ is complete to $Y$, and thus $v$ is nonadjacent to some $z_j$. Note that $y_j$ exists since $s =\lceil k/2\rceil\ge \lfloor k/2\rfloor=t$. Thus, $\widetilde G[X\cup\{y_j,z_j,v,w\}]$ is the fifth graph, which implies that $G$ contains the fifth graph. \hfill $\Box$\medskip

Theorem \ref{thm:remedy} simplifies if $G$ or $\ov G$ is claw-free; we omit the straightforward proof.

\begin{corollary}\label{cor:remedy}
Let $G$ be doubly connected such that either $G$ or $\ov G$ is claw-free. Let $H\subset G$ with $0<||H||< \binom{|H|}{2}$. If some vertex of $G-H$ is complete to $H$, then there exists an $H$-remedy $J\subseteq G$.
\end{corollary}

Our final result in this section addresses irreducible claw-free graphs $G$. For any connected induced subgraph $H$ with a cut vertex $s$ (so $H$ is a simplicial sum and hence reducible), we determine all augmentations of $H$ within $G$ that ensure $s$ is no longer a cut vertex.

\begin{theorem}\label{thm:sumblock}
Let $G$ be claw-free and irreducible. Let $H\subset G$ be the simplicial sum of $H_1, H_2$ over their common simplicial vertex $s$. For each $u\in U=V(G-H)$, let $Z_u= N_G(u)\cap V(H)$, and for each $i=1,2$, let $N_i=N_G(s)\cap V(H_i)$, $R_i=V(H_i)-N_{H_i}[s]$. If $R_1\ne\emptyset$, then at least one of the following must hold. 
\begin{enumerate}[label=(\arabic*)]
\parskip -4pt
\item $\exists\ u\in U$ such that $Z_u\cap V(H_1-s)\ne\emptyset\ne N_1-Z_u$ and $s\in Z_u$;
\item $\exists\ u\in U$ such that $Z_u\cap R_1\ne\emptyset\ne Z_u\cap V(H_2-s)$ and $s\in Z_u$;
\item $\exists\ u\in U$ such that $Z_u\cap V(H_1-s)\ne\emptyset\ne Z_u\cap V(H_2-s)$ and $s\not\in Z_u$;
\item $\exists\ u, v \in U$ such that $Z_u\cap V(H_1-s)\ne \emptyset\ne Z_v\cap V(H_1-s)$, $s\in Z_u\cap Z_v$, and $uv\not\in E(G)$;
\item $\exists\ u, v \in U$ such that $Z_u\cap R_1\ne\emptyset=Z_v\cap V(H_1-s)$, $s\in Z_u$, and $uv \in E(G)$;
\item $\exists\ u, v \in U$ such that $Z_u\cap V(H_1-s)\ne\emptyset=Z_v\cap V(H_1-s)$, $s\not\in Z_u$, and $uv \in E(G)$;
\item $\exists\ u, v \in U$ such that $Z_u\cap V(H_1-s)\ne\emptyset = Z_v\cap  R_1$, $s\in Z_v-Z_u$, and $uv \in E(G)$.
\end{enumerate} 
\end{theorem}

\noindent{\bf Proof.} Suppose none of (1)-(7) holds. We will prove that $G$ is reducible, which yields a desired contradiction. We begin by defining a partition of $V(G-H)$: \\ 
\indent $A_0=\{u\in V(G-H): Z_u\cap V(H_1-s)=\emptyset\}$,\\ 
\indent $A_1=\{u\in V(G-H-A_0): us\in E(G)$ and $Z_u\cap R_1=\emptyset\}$,\\ 
\indent $A_2=\{u\in V(G-H-A_0): us\in E(G)$ and $Z_u\cap R_1\ne \emptyset\}$, and \\ 
\indent $A_3=V(G-H)-A_0-A_1-A_2$. \\ 
Next, we establish a few properties of this partition.

\noindent {\bf Claim 1.} {\it $N_1\cup A_1\cup A_2\cup \{s\}$ is a clique.}\\ 
\indent First, the definitions of $N_1,A_1,A_2$ imply that $s$ is complete to  $N_1\cup A_1\cup A_2$. Next, the definitions of $A_0,A_1,A_2,A_3$ imply that every $u\in A_1\cup A_2\cup A_3$ satisfies $Z_u\cap V(H_1-s)\ne\emptyset$. As such, if some $u\in A_1\cup A_2$ is not complete to $N_1$, then $u$ satisfies (1); if some distinct $u,v\in A_1\cup A_2$ are nonadjacent, then the pair $u,v$ satisfies (4). Thus the claim holds. 

\noindent {\bf Claim 2.} {\it $N_G(R_1\cup A_3)\subseteq N_1\cup A_2$.}\\ 
\indent To see this, first observe that by the definitions of $R_1$, $A_0$, and $A_1$, every vertex of $R_1$ is anticomplete to $V(H_2)\cup A_0\cup A_1$; hence $N_G(R_1)\subseteq N_1\cup A_2\cup A_3$. Next, we consider vertices $u\in A_3$. By the definition of $A_3$, $s\not\in Z_u$. If $u$ is adjacent to a vertex of $H_2$, then $u$ satisfies (3); if $u$ is adjacent to some $v\in A_0\cup A_1$, then the pair $u,v$ satisfies (6) or (7). Consequently, $N_G(A_3)\subseteq N_1\cup A_2\cup R_1$ and thus the claim holds.

\noindent {\bf Claim 3.} {\it $N_G(V(H_2-s)\cup A_0)\subseteq A_1\cup\{s\}$.}\\ 
\indent The definition of $N_1$ implies that $N_1$ is anticomplete to $H_2-s$. Meanwhile, if some $u\in A_2$ has a neighbor in $H_2-s$, then $u$ satisfies (2). Together with Claim 2, we conclude that $N_G(V(H_2-s))\subseteq A_0\cup A_1\cup\{s\}$. Next, we consider vertices $v\in A_0$. The definition of $A_0$ implies that $v$ is anticomplete to $N_1$. Moreover, if $v$ is adjacent to some $u\in A_2$, then the pair $u,v$ satisfies (5). Again, together with Claim 2, we conclude that $N_G(A_0)\subseteq V(H_2)\cup A_1$, and thus the claim holds.

Note that $B_0=A_0\cup V(H_2-s)$, $B_1=A_1\cup\{s\}$, $B_2 = A_2\cup N_1$, and $B_3=A_3\cup R_1$ form a partition of $V(G)$. From the definition of simplicial sum and the assumption of the theorem on $R_1$, we see that $N_1$, $N_2$, $R_1$ are non-empty, and thus $B_0,B_1,B_2,B_3$ are non-empty. From Claims 1-3 we see that $B_1\cup B_2$ is a clique, $N_G(B_0)\subseteq B_1$, and $N_G(B_3)\subseteq B_2$. Therefore, $G$ is either a simplicial sum (if $|B_1|=1$ or $|B_2|=1$, since $G$ is claw-free) or a simplicial join (if $|B_1|$, $|B_2|\ge2$). This contradiction proves the theorem. \hfill $\Box$

\section{Splitter theorems}

For any graph $H$, an $H$-{\it add} is a graph obtained from $H$ by adding a new vertex $v$ adjacent to at least one vertex of $H$; we call $v$ the {\it added vertex}. Similarly, an $H$-{\it add$_1$} is an $H$-add where the added vertex is mixed on $H$. Note that if $H$ is connected, then every $H$-add is connected. Likewise, if $H$ is doubly connected, then every $H$-add$_1$ is doubly connected; moreover, these are the only $H$-adds that are doubly connected. 

If $H$ is a proper induced subgraph of $G$, where both $H$ and $G$ are connected, then $G$ clearly contains an induced $H$-add. Consequently,  $G$ can be constructed from $H$ by iteratively applying the add operation while keeping all intermediate graphs connected. 
However, the analogous statement fails for double connectivity, as add$_1$ alone is insufficient to generate $G$ from $H$. In particular, a $P_4$-augmentation of $H$ cannot be generated from $H$ via add$_1$ operations. It turns out that the operation of $P_4$-augmentation is also required, and that these two operations together are sufficient, as stated in Theorem \ref{thm:chain-dc}.

The objective of this section is to establish analogous results for prime graphs. In particular, we show that every connected prime graph can be constructed from any nontrivial connected prime induced subgraph by add operations that preserve connectedness and primeness. We also prove a parallel result for doubly connected graphs. Finally, we introduce another variant of the add operation that yields a similar result for irreducible prime graphs in a key special case.

\subsection{Connected graphs}

Our goal in this subsection is to prove the following two results. 

\begin{theorem}\label{thm:chain-prime}
Let $G$ and $H$ be prime graphs with $H \prec G$. Then there exists a sequence $H \cong G_0\prec G_1\prec \cdots \prec G_k\cong G$ of prime graphs such that $|G_i| = |G_{i-1}|+1$ for all $i\in [k]$.
\end{theorem}

\begin{theorem}\label{thm:chain-connp}
Let $G$ and $H$ be connected prime graphs with $H \prec G$ and $|H|>1$. Then there exists a sequence $H \cong G_0\prec G_1\prec \cdots \prec G_k\cong G$ of connected prime graphs such that $|G_i| = |G_{i-1}|+1$ for all $i\in [k]$.
\end{theorem}

\begin{figure}[ht]
\centering
\includegraphics[scale=0.72]{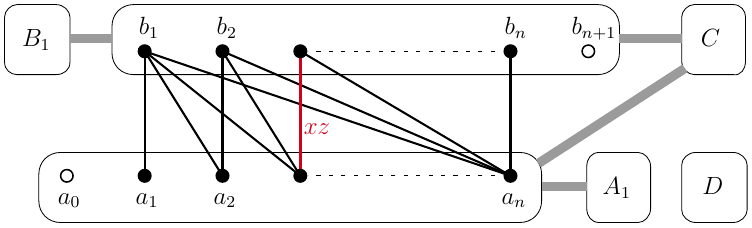}
\captionsetup{justification=centering}
\caption{\small Structure of $G$; thick lines indicate pairs of sets that are complete to each other.}
\label{f1}
\end{figure}

A key ingredient of our proof is a {\it quadruple} $(G,x,y,z)$, where $G$ is a graph and $x,y,z\in V(G)$ satisfy the following conditions: $x$ and $y$ are twins in $G-z$, $zxy$ is an induced path, and this path is not a connected component of $G$. For such a quadruple, we prove that $V(G)$ admits a partition $(A_0,B_0,A_1,B_1,C,D)$ satisfying the following properties (see Figure \ref{f1}):
\begin{enumerate}
\parskip -3pt
\item[(P1)] There exists an integer $n\ge 1$ such that $A_0=\{a_0, a_1, \ldots, a_n\}$ and  $B_0=\{b_1, \ldots, b_n,b_{n+1}\}$, where $a_0$ and $b_{n+1}$ may be absent. Furthermore, the edge $xz$ equals $a_kb_k$ for some $k\in [n]$.
\item[(P2)] $N[a_i] = \{b_1,\ldots, b_i\}\cup  A_0\cup A_1 \cup C$, for every $a_i\in A_0$, and \\ 
$N[b_i] = \{a_i,\ldots, a_n\}\cup B_0\cup B_1\cup C$, for every $b_i\in B_0$.
\end{enumerate}
Note that both $A_0$ and $B_0$ are cliques. In addition, all edges incident to a vertex of $A_0\cup B_0$ are specified by (P2); the only other edges of $G$ are between vertices in $A_1\cup B_1\cup C\cup D$. 

\begin{lemma}\label{lem:partition}
Let $(G,x,y,z)$ be a quadruple such that $G_0=G-x-z$ is prime. Then $V(G)$ admits a partition $(A_0,B_0,A_1,B_1,C,D)$ satisfying {\rm (P1)} and {\rm (P2)}. Moreover, if the partition is chosen to maximize $|A_0\cup B_0|$, then $|A_0\cup B_0|\ge 3$, and the following hold.
\begin{enumerate}[label=(\arabic*)]
\parskip -4pt
\item For each $i\in [n]$, $G-a_i-b_i\cong G_0$. Moreover, if $a_0$ is present, then $G-a_0-b_1\cong G_0$; if $b_{n+1}$ is present, then $G-a_n-b_{n+1}\cong G_0$. 
\item If both $a_0$ and $b_{n+1}$ are present, then $G-a_0$ is a prime $(G-a_0-b_1)$-add, and $G-b_{n+1}$ is a prime $(G-a_n-b_{n+1})$-add.  
\item If $a_0$ is absent, then $G-b_1$ is a prime $(G-a_1-b_1)$-add$_1$, unless $d_G(a_1)=1$, in which case $n=1$, $b_{n+1}$ is present, and $G-b_{n+1}$ is a prime $(G-a_n-b_{n+1})$-add$_1$. Symmetrically, if $b_{n+1}$ is absent, then $G-a_n$ is a prime $(G-a_n-b_n)$-add$_1$, unless $d_G(b_n)=1$, in which case $n=1$, $a_0$ is present, and $G-a_0$ is a prime $(G-a_0-b_1)$-add$_1$. 
\end{enumerate}
\end{lemma}

\noindent{\bf Proof.} We first prove that $V(G)$ admits a partition satisfying (P1), (P2), and $|A_0\cup B_0|\ge3$. Let $A_0=\{a_0,a_1\}$ and $B_0=\{b_1\}$, where $a_0=y$, $a_1=x$, and $b_1=z$. We then partition the remaining vertices of $G$ according to their adjacency to $y$ and $z$: $A_1= N(y)-N(z)$, $B_1=N(z)-N(y)$, $C=N(y)\cap N(z) - x$, and $D=V(G)-N[y]-N[z]$. By construction, $|A_0\cup B_0|\ge3$, and (P1) holds with $n=1$. In addition, 
\begin{align*}
N[b_1] = N[z] &= \{z\} \cup (N(z)-N(y))\cup (N(y) \cap N(z)-x) \cup \{x\} = \{a_1\} \cup B_0 \cup B_1 \cup  C, \\ 
N[a_0] = N[y] &= \{y\} \cup (N(y)-N(z)) \cup (N(y) \cap N(z)-x) \cup \{x\} = A_0 \cup A_1 \cup C, \text{ and}\\ 
N[a_1] = N[x] &= \{z\} \cup N[y] = \{b_1\} \cup A_0 \cup A_1\cup C. 
\end{align*} 
Thus (P2) also holds, establishing the existence of the desired partition. Among all such partitions, choose one that maximizes $|A_0\cup B_0|$. It remains to show that this partition satisfies (1)-(3).

To prove (1), we first verify that $G-a_i-b_i\cong G-a_{i+1}-b_{i+1}$ for every $i\in [n-1]$. Since $a_ib_i, a_{i+1}b_{i+1}\in E(G)$, $N[a_i] = N[a_{i+1}] -b_{i+1}$, and $N[b_{i+1}] = N[b_i]-a_i$, the map that swaps $a_i$ with $a_{i+1}$ and $b_i$ with $b_{i+1}$ while fixing all other vertices is an isomorphism between $G-a_i-b_i$ and $G-a_{i+1}-b_{i+1}$. Since $G_0$ is one of these graphs, it follows that $G-a_i-b_i\cong G_0$ for every $i\in [n]$, as required. Now suppose that $a_0$ is present, then $a_0$ and $a_1$ are twins in $G-b_1$, and hence $G-a_0-b_1\cong G-a_1-b_1\cong G_0$. By symmetry, if $b_{n+1}$ is present, then $G-a_n-b_{n+1}\cong G_0$. This completes the proof of (1).

To prove (2), by symmetry it suffices to prove the second assertion. Let $H=G-a_n-b_{n+1}$ and $J=G-b_{n+1}$. Then $J$ is an $H$-add with added vertex $a_n$, as $a_n$ is adjacent to $b_n\in V(H)$. It remains to show that $J$ is prime. Suppose otherwise. Since $H\cong G_0$ is prime,  by Observation~\ref{obs:basic}(2), $a_n$ has a twin $w$ in $J$. We locate $w$ as follows. Since $a_0$ is adjacent to $a_n$ and anticomplete to $B_0$, we have $w\notin B_0$. Since $b_n$ is adjacent to $a_n$ and anticomplete to $A_0\cup A_1-a_n$, we have $w\notin A_0\cup A_1$. Finally, $w\notin B_1\cup D$, as $a_n$ is anticomplete to this set. Hence $w\in C$. Since $N_G[w]=N_J[w]\cup\{b_{n+1}\}=$ $N_J[a_n]\cup\{b_{n+1}\} = N_G[a_n] \cup\{b_{n+1}\}=B_0\cup A_0\cup A_1\cup C$, moving $w$ from $C$ to $A_0$ and setting $a_{n+1}=w$ yields a new partition satisfying (P1) and (P2). This increases $|A_0\cup B_0|$, contradicting the maximality of the chosen partition. Therefore, $J$ is prime, completing the proof of (2).  
 
Finally we prove (3). Again by symmetry, it suffices to prove the first assertion. We first consider the case $d_G(a_1)=1$. Then $A_0=\{a_1\}$, so $n=1$ and $B_0\subseteq \{ b_n, b_{n+1}\}$. Since $|A_0\cup B_0|\ge3$, the maximality of $|A_0\cup B_0|$ implies that $b_{n+1}$ is present. Because $(G,x,y,z)$ is a quadruple, $x$ and $y$ are twins in $G-z$, and $zxy$ is an induced path that is not a component of $G$. Thus, from $a_1b_1=xz$ and $d_G(x)\ge2$, we obtain $a_1=z$ and $b_1=x$. Since $d_G(a_1)=1$, we also have $A_1=C=\emptyset$, so $b_1$ and $b_2$ are twins in $G-a_1$. As $x$ and $y$ are twins in $G-z$, it follows that $b_2=y$, for otherwise $b_2$ and $y$ are twins in $G_0$, contradicting our hypothesis that $G_0$ is prime. Consequently, $B_1\ne\emptyset$, since the path $zxy$ is not a component of $G$. Let $H=G-a_n-b_{n+1}$ and $J=G-b_{n+1}$. Then $J$ is an $H$-add$_1$, as $a_n$ is adjacent to $b_n\in V(H)$ and anticomplete to $B_1$. Moreover, $J$ is prime by Observation \ref{obs:basic}(2), since $z$ and its only neighbor $x$ are not twins in $J$ (because $B_1\ne\emptyset$). Hence $J$ is a prime $H$-add$_1$, as required. 

It remains to consider the case $d_G(a_1)\ge 2$. Let $H=G-a_1-b_1$ and $J=G-b_1$. Since $d_G(a_1)\ge 2$, we have $d_J(a_1)\ge 1$, so $J$ is an $H$-add with added vertex $a_1$. Moreover, $b_2$ exists because $|B_0|\ge \lceil |A_0\cup B_0|/2 \rceil \ge \lceil 3/2\rceil =2$. Thus $J$ is an $H$-add$_1$ since $a_1$ is nonadjacent to $b_2$. It remains to show that $J$ is prime. Suppose not. Since $H\cong G_0$ is prime, by Observation \ref{obs:basic}(2), $a_1$ has a twin $w$ in $J$. Note that $w\ne a_i$ for any $i\ge2$ since $b_i$ is adjacent to $a_i$ but not to $a_1$. Also, $w\notin C$, for otherwise $b_2$ would be adjacent to $w$ but not to $a_1$. Hence $w\in A_1$. Since $N_G[w]=$ $N_J[w]=N_J[a_1]=N_G[a_1]-b_1=A_0\cup A_1\cup C$, moving $w$ from $A_1$ to $A_0$ and setting $a_0=w$ yields a new partition satisfying (P1) and (P2). This increases $|A_0\cup B_0|$, contradicting the maximality of the chosen partition. Therefore, $J$ is a prime $H$-add$_1$, completing the proof of the lemma. \hfill $\Box$\medskip

We prove a stronger result below, which implies Theorem \ref{thm:chain-connp} immediately.

\begin{lemma}\label{lem:primeadd}
Let $G$ be connected and prime. Let $H$ be nontrivial and prime with $H\prec G$. Then there exists a prime $H$-add $J$ such that $J\preceq G$.
\end{lemma}

\noindent{\bf Proof}. Since $H \prec G$, $G$ has an induced subgraph $G_0\cong H$. As $G$ is connected, it contains a $G_0$-add $G_0^+$, with added vertex $x$. If $G_0^+$ is prime, then $J=G_0^+$ satisfies the conclusion. Hence, we assume that $G_0^+$ is not prime. Since $G_0$ is prime, by Observation \ref{obs:basic}(2), $x$ has a twin $y$ in $G_0^+$. As $G$ is prime, it has a vertex $z$ adjacent to exactly one of $x$ and $y$, say $x$. Let $G_0^* = G_0+\{x,z\}$. Then $zxy$ is an induced path in $G_0^*$. If $zxy$ is not a component of $G_0^*$, then $(G_0^*, x, y, z)$ satisfies the hypotheses of Lemma \ref{lem:partition}. By this lemma, there exists a prime graph $J\subseteq G_0^*\subseteq G$ such that $J$ is isomorphic to an $H$-add. Hence, $J$ satisfies the conclusion. 

It remains to consider the case in which $zxy$ is a component of $G_0^*$. In this case $y$ is an isolated vertex of $G_0$. Let $F=G_0-y$. Then $|F|=|G_0|-1=|H|-1>0$. Since $G$ is connected, it has a shortest path $v_1v_2\cdots v_k$ from $F$ to $y$, where $k\ge3$ and $v_1\in V(F)$. Let $G_1=F+v_3$. Then $G_1\cong G_0\cong H$ and $G_1^+= G_1 +v_2$ is a $G_1$-add since the added vertex $v_2$ is adjacent to $v_3$. In addition, $v_2$ has no twins in $G_1^+$ since $G_1^+[N_{G_1^+}(v_2)]$ is disconnected. Hence $G_1^+$ is prime, and $J=G_1^+$ satisfies the conclusion. This completes the proof. \hfill $\Box$\medskip

\noindent{\bf Proof of Theorem \ref{thm:chain-prime}}. It suffices to show that there exists a prime graph $J$ satisfying $H\prec J\preceq G$ and $|J|=|H|+1$, since repeated application of this assertion yields the desired sequence. Choose $G_0\subseteq G$ with $G_0\cong H$ such that, whenever a nontrivial component $G'$ of $G$ contains exactly one vertex $v_0$ of $G_0$, the vertex $v_0$ has a non-neighbor in $G'$. This is possible because $G'$ is prime and hence non-complete, while replacing $v_0$ by any other vertex of $G'$ preserves the condition $G_0\cong H$.

If some $v\in V(G-G_0)$ is anticomplete to $G_0$, then $J=G_0 + v$ is obtained from the prime graph $G_0$ by adding an isolated vertex. It follows that $J$ is prime, and therefore is a desired graph. Thus we may assume that every vertex of $G-G_0$ has a neighbor in $G_0$. 

Let $G'$ be a component of $G$ meeting $G-G_0$, and let $G_0'=G_0\cap G'$. By the preceding assumption, $G_0'$ is not empty. If $G_0'$ consists of only one vertex $v_0$, then the choice of $G_0$ implies that $v_0$ has a non-neighbor $v$ in $G'$. Note that $v$ is anticomplete to $G_0$, contradicting our assumption. Hence $G_0'$ is nontrivial. In addition, $G_0'$ is prime because $G_0$ is prime and $G_0'\subseteq G_0$ is anticomplete to $G_0-G_0'$. 
Applying Lemma \ref{lem:primeadd} to $G_0'$ and $G'$, we obtain a prime $G_0'$-add $J'\preceq G'$. Then $J = J' \cup (G_0-G_0')$ is a desired graph, completing the proof. 
\hfill $\Box$

\subsection{Doubly connected graphs}

As noted in the beginning of this section, every doubly connected graph can be generated from any of its doubly connected induced subgraphs via add$_1$ and $P_4$-augmentation operations.

\begin{theorem}\label{thm:chain-dc}
Let $G$ and $H$ be doubly connected with $H \prec G$. Then there exists a sequence $H \cong G_0\prec G_1\prec \cdots \prec G_k\cong G$ of doubly connected graphs such that each $G_i$ $(i=1,...,k)$ is either a $G_{i-1}$-add$_1$ or a $P_4$-augmentation of $G_{i-1}$.
\end{theorem}

We omit the proof of this result since it is an immediate consequence of Theorem \ref{thm:remedy}. We remark that, by Corollary \ref{cor:remedy}, if either $G$ or $\ov G$ is claw-free, then $P_4$-augmentation can be eliminated from the statement.

The goal of this subsection is to establish a doubly connected analogue of Theorem \ref{thm:chain-connp}. However, the corresponding result is considerably more involved for two fundamental reasons. First, unlike the connected case, the add$_1$ operation alone is insufficient; $P_4$-augmentations are also required. Second, a new obstacle arises in the doubly connected setting, as we explain below.
 
 \begin{figure}[ht]
\centering
\includegraphics[scale=0.68]{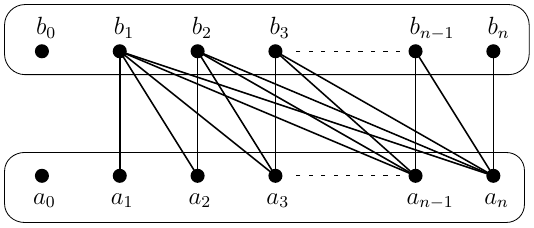}
\caption{The exceptional graphs $T_n$.}
\label{fig:Tn}
\end{figure}

Let $\mathcal T=\{T_n: n\ge1\}$, where $T_n$ is the graph (illustrated in Figure \ref{fig:Tn}) consisting of two disjoint cliques $\{a_0, a_1, . . ., a_n\}$ and $\{b_0, b_1, . . . , b_n\}$ such that edges between them are exactly those $a_ib_j$ satisfying $1\le j\le i\le n$. Note that $T_1$ is $P_4$ and $T_2$ is $\mathbb L_3$ (See Figure \ref{fig:line7}). The two cliques are symmetric and each $T_n$ is both prime and doubly connected. On the other hand, $\ov{T_n}-b_0$ contains the isolated vertex $a_n$, while $T_n-b_i$ ($i\ge1$) contains the twins $a_{i-1}$ and $a_i$. Therefore, if $n\ge 2$, no graph strictly between $T_{n-1}$ and $T_n$ can be simultaneously prime and doubly connected. As a result, $T_n$ cannot be obtained from $T_{n-1}$ via add$_1$ operations and $P_4$-augmentations while preserving these two properties. Fortunately, the family $\cal T$ is the only exception.

\begin{theorem}\label{thm:primedc}
Let $G$ and $H$ be doubly connected prime graphs with $G\not\in\cal T$ and $H\prec G$. Then there exists a sequence $H \cong G_0\prec G_1\prec \cdots \prec G_k\cong G$ of doubly connected prime graphs such that each $G_i$ $(i=1,...,k)$ is either a $G_{i-1}$-add$_1$ or a $P_4$-augmentation of $G_{i-1}$.
\end{theorem}

As in Theorem \ref{thm:chain-dc}, if either $G$ or $\ov G$ is claw-free, then the $P_4$-augmentation operation can be eliminated from Theorem \ref{thm:primedc}. The proof of Theorem \ref{thm:primedc} occupies the remainder of this subsection. Since one part of the argument is used repeatedly, we isolate it as a lemma.

\begin{lemma}\label{lem:primedc}
Let $G$ and $H$ be doubly connected prime graphs with $H\prec G$. Then there exists a doubly connected prime graph $J\preceq G$ that is either an $H$-add$_1$ or a $P_4$-augmentation of $H$, unless $H=T_t$ and $T_{t+1}\preceq G$ for some $t\ge 1$. 
\end{lemma}

\noindent{\bf Proof}. By Theorem \ref{thm:chain-dc}, there exists a doubly connected graph $H^+\preceq G$ that is either an $H$-add$_1$ or a $P_4$-augmentation of $H$. Since a routine verification shows that every $P_4$-augmentation of a prime graph is again prime, $J=H^+$ is a desired graph unless $H^+$ is a non-prime $H$-add$_1$, which we now assume.

Without loss of generality, let $G_0\subseteq G_0^+\subseteq G$, where $G_0\cong H$ and $G_0^+\cong H^+$. Let $x$ be the added vertex in obtaining $G_0^+$. Since $G_0$ is prime whereas $G_0^+$ is not, Observation \ref{obs:basic}(2) implies that $x$ has a twin $y$ in $G_0^+$. As $G$ is prime, it has a vertex $z$ adjacent to exactly one of $x$ and $y$, say $x$. Let $G_0^*=G_0 + \{x,z\}$. Since $G_0$ is doubly connected, the induced path $zxy$ is not a component of $G_0^*$. Hence Lemma \ref{lem:partition} applies to the quadruple $(G_0^*, x, y, z)$. 

In the remainder of this proof we use the notation of Lemma \ref{lem:partition} (see Figure \ref{f1}). Suppose first that either $a_0$ or $b_{n+1}$ is absent. By Lemma \ref{lem:partition}(3), there are four scenarios and in each case some $J\subseteq G_0^*$ is a prime $G'$-add$_1$ for some $G'\subseteq G_0^*$. By Lemma \ref{lem:partition}(1), $G'\cong G_0\cong H$ in each case, so $J$ satisfies the requirement in each case. We may therefore assume that both $a_0$ and $b_{n+1}$ are present. By Lemma \ref{lem:partition}(2), $G_0^*-a_0$ is a prime $(G_0^*-a_0-b_1)$-add, and $G_0^*-b_{n+1}$ is a prime $(G_0^*-a_n-b_{n+1})$-add. If $A_1\cup B_1\cup D\ne \emptyset$, then either $J=G_0^*-a_0$ or $J=G_0^*-b_{n+1}$ satisfies the requirement, since either $b_1$ has a non-neighbor in $G_0^*-a_0$ or $a_n$ has a non-neighbor in $G_0^*-b_{n+1}$. Hence we may further assume that $A_1\cup B_1\cup D= \emptyset$. We claim that $C=\emptyset$ as well. Indeed, if $C\ne\emptyset$, then since $C$ is complete to $A_0\cup B_0$, the complement of $G_0^*-x-z\cong H$ is disconnected, contradicting the assumption that $H$ is doubly connected. Consequently, $G_0^*\cong T_n$ and $H\cong G_0^*-x-z\cong T_{n-1}$, which establishes the exceptional case and completes the proof.\hfill $\Box$\medskip

\noindent{\bf Proof of Theorem \ref{thm:primedc}}. It suffices to find a doubly connected prime graph $J\subseteq G$ that is isomorphic to either an $H$-add$_1$ or a $P_4$-augmentation of $H$, because repeatedly applying this result yields a chain required in the Theorem. To find a desired $J$, by Lemma \ref{lem:primedc}, we only need to consider the case that $H=T_t$ for some $t\ge 1$. Choose $T_n\preceq G$ such that $n$ is maximum. Then $n\ge t$. Since $G\notin\cal T$, we have $T_n\prec G$. By applying Lemma \ref{lem:primedc} to $T_n$ and $G$ we conclude that there exist $G_0\subseteq G_1\subseteq G_1^+\subseteq G$ such that $G_0\cong T_t$, $G_1\cong T_n$, and $G_1^+$ is either a prime $G_1$-add$_1$ or a $P_4$-augmentation of $G_1$. Since every $P_4$-augmentation of $T_n$ contains a $P_4$-augmentation of $T_t$, we may assume that $G_1^+$ is a prime $G_1$-add$_1$ with added vertex $v$. If $n=t$ then $G_0=G_1$ and $J=G_1^+$ is a desired graph. So we assume $n>t$. In particular, $n\ge2$. Let the vertices of $G_1$ be $a_0, b_0, ..., a_n,b_n$ as in the definition of $T_n$. We prove that either $G_1^+$ contains an induced subgraph isomorphic to a prime $G_0$-add$_1$ or for some  $i\in [n]$, $G_1^+-a_i-b_i$ is a prime $(G_1-a_i-b_i)$-add$_1$. Note that repeatedly applying this result leads to a desired $J$ and thus proves the theorem. 

Suppose $v$ is adjacent to both $a_0$ and $b_0$. Since $v$ is mixed on $G_1$, we may choose a non-neighbor $x$ of $v$. Let $j$ be the index with $x\in\{a_j,b_j\}$. Since $n\ge 2$, there exists $i\in [n]$ with $i\ne j$. Then $G_1^+-a_i-b_i$ is a $(G_1-a_i-b_i)$-add$_1$. In addition, $v$ has no twin in $G_1^+-a_i-b_i$ since $v$ is adjacent to both $a_0,b_0$. So $G_1^+-a_i-b_i$ is also prime, implying that $i$ satisfies the requirement. The argument is similar for the case that $v$ is nonadjacent to both $a_0$ and $b_0$. Let $x$ be a neighbor of $v$ and let $x\in\{a_j,b_j\}$. Then any $i\in [n]-\{j\}$ satisfies the requirement. 

It remains to consider the case that $v$ is adjacent to exactly one of $a_0,b_0$, say $b_0$. If some $b_j$ is a non-neighbor of $v$, then any $i\in [n]-\{j\}$ satisfies the requirement. So $v$ is adjacent to every $b_j$. Since $G_1^+$ is prime, $v$ is not a twin of $b_0$, so some $a_j$ is a neighbor of $v$. Let $j$ be the smallest index with this property. Since $v$ and $b_j$ are not twins, some $a_k$ is a non-neighbor of $v$ with $k>j$. If there exists $i\in [n]-\{j,k\}$, then $i$ satisfies the requirement. If no such $i$ exists, then $n=2$, which means $t=1$. In this case, $G_0\cong T_t=T_1\cong G_1[b_0,b_2,a_2,a_1]$, and $G_1^+[v,b_0,b_2,a_2,a_1]$ is a prime $G_1[b_0,b_2,a_2,a_1]$-add$_1$. This completes the proof. \hfill $\Box$

\subsection{Irreducible graphs}

All graphs considered in this subsection are claw-free. 
The previous two subsections established splitter theorems for connected prime graphs and doubly connected prime graphs (Theorems \ref{thm:chain-connp} and \ref{thm:primedc}). A natural question is whether analogous results hold for irreducible graphs, either in a prime setting or, in general. As a first step, we ask whether every irreducible graph $G$ containing an irreducible induced subgraph $H$ must contain an irreducible local augmentation of $H$. 

Unfortunately, the answer is negative in both settings. Indeed, let $G$ be the union of $H$ and a path $P$ whose ends are simplicial vertices of $H$ and whose internal vertices are disjoint from $H$. If $|P|\ge4$, then every graph strictly between $H$ and $G$ is reducible. Thus no direct analogue of the previous splitter theorems can hold in full generality.

This example motivates us to restrict our attention to irreducible graphs with at most one simplicial vertex. We prove that every claw-free irreducible graph $G$ containing an irreducible induced subgraph $H$ with at most one simplicial vertex must contain an irreducible $H$-add. We also establish the corresponding result for prime graphs.
These theorems, however, are not themselves splitter theorems, since the resulting $H$-add may introduce additional simplicial vertices, preventing the construction from being iterated. To make the construction iterative, we introduce at the end of this subsection a modified add operation that yields a chain version of the above results.
 
We begin with a simple observation that identifies the only two scenarios in which an $H$-add is reducible for an irreducible graph $H$.

\begin{lemma}\label{lem:irradd}
Let $H\ne K_2$ be irreducible and let $G$ be a claw-free $H$-add with added vertex $v$. Then $G$ is reducible if and only if $H$ has a simplicial vertex $s$ such that either $v$ is adjacent only to $s$ or $v$ is a twin of $s$.
\end{lemma}

\noindent{\bf Proof.} For the backward implication, we first note from Observation \ref{obs:basic}(3) that $d_H(s)>1$. As such, if $v$ is adjacent only to $s$, then $G$ is a simplicial sum of $H$ and $K_2$, and if $v$ is a twin of $s$, then $G$ is a simplicial join of $H$ and $K_3$. 
For the forward implication, note that $G$ is connected and thus is a simplicial-sum or -join. Suppose first that $G$ is a simplicial sum of $G_1$ and $G_2$ over their common simplicial vertex $s$. Without loss of generality, let $v\in V(G_2)$. If $|G_2|\ge3$, then either $v=s$ and $G-v$ is disconnected, or $v\ne s$ and $G-v$ is a simplicial sum of $G_1$ and $G_2-v$. In either case, $G-v$ is reducible, contradicting the irreducibility of $H=G-v$. Hence $|G_2|=2$ and $v$ is adjacent only to $s$. Now suppose that $G$ is a simplicial join of $G_1,G_2$ over their respective simplicial vertices $v_1,v_2$, and let $v\in V(G_2)$. If $|N_{G_2}(v_2)-v|\ge2$, then $G-v$ is a simplicial join of $G_1$ and $G_2-v$. If $|N_{G_2}(v_2)-v|=1$, say $N_{G_2}(v_2)=\{v,s\}$, and $V(G_2)-N_{G_2}[v_2]\ne\emptyset$, then $G-v$ is a simplicial sum with cut vertex $s$. Both cases contradict the irreducibility of $H$. Hence $|N_{G_2}(v_2)|=2$ and $V(G_2)=N_{G_2}[v_2]$, so $v$ is a twin of the simplicial vertex $s$ of $H$. \hfill $\Box$\medskip

Our first theorem is stated in a slightly stronger form. Together with the preceding lemma, it clearly implies that every claw-free irreducible graph $G$ properly containing an irreducible graph $H$ with at most one simplicial vertex must contain an irreducible $H$-add.
 
\begin{theorem}\label{thm:addirr}
Let $H\not\in \{K_1,K_2\}$ be a graph with a specified simplicial vertex $s$, and let $H\prec G$ for a claw-free irreducible graph $G$. Then there exists an $H$-add $H^+$ with added vertex $v$ such that $H^+\preceq G$ and $N_{H^+}(v)\not\in\{\{s\},N_H[s]\}$.
\end{theorem}

\noindent{\bf Proof}. It suffices to find $J\subseteq G$ isomorphic to the graph $H^+$ stated in the theorem. Let $G_0\subseteq G$ with $G_0\cong H$, and let $x\in V(G_0)$ correspond to $s$. Let $R=\{r\in V(G-G_0): N_G(r)\cap V(G_0)\ne\emptyset\}$, and define $R_1=\{r\in R: N_G(r)\cap V(G_0)=\{x\}\}$ and $R_2=\{r\in R: N_G(r)\cap V(G_0)=N_{G_0}[x]\}$. If there exists $v\in R-R_1-R_2$, then $J=G_0+v$ satisfies the conclusion. Hence we may assume that $R=R_1\cup R_2$. Since every vertex of $R$ is adjacent to $x$, any component of $G_0$ not containing $x$ would also be a component of $G$, contradicting the irreducibility of $G$. Thus $G_0$ is connected. Since $H\notin \{K_1,K_2\}$, we have $|G_0|\ge 3$, so $N_{G_0}[x]\ne\{x\}$, and thus $R_1\cap R_2=\emptyset$. 

We next show that $R_2\ne\emptyset$. Otherwise $R=R_1$, so $x$ separates $G_0-x$ from $G-G_0$. Since $|G|>|G_0|>1$, $x$ is a cut vertex of $G$. By Observation \ref{obs:basic}(3), $G$ is a simplicial sum, contradicting its irreducibility. We also claim that $d_{G_0}(x)\ge2$. Indeed, if $N_{G_0}(x)$ consists of a single vertex $u$, then, as $|G_0|\ge3$ and $R=R_1\cup R_2$, the vertex $u$ separates $G_0-\{x,u\}$ from $G-(G_0-x)$, making $u$ a cut vertex of $G$. Again, $G$ is a simplicial sum, a contradiction.  

Suppose that $R_2$ is a clique. Consider the partition $(X_1,X_2,X_3,X_4) = (V(G_0)-N_{G_0}[x], N_{G_0}(x)$, $\{x\}\cup R_2$, $V(G-G_0)-R_2)$ of $V(G)$. Then $X_2\cup X_3$ is a clique, $X_1$ is anticomplete to $X_3\cup X_4$, $X_4$ is anticomplete to $X_1\cup X_2$, and $\min\{|X_2|,|X_3|\}\ge2$. Hence $G$ is a simplicial join, contradicting its irreducibility. Therefore $R_2$ is not a clique. Choose nonadjacent vertices $u,v\in R_2$ and let $G_1=$ $(G_0-x)+u$. Since $N_G(u)\cap V(G_0)=N_{G_0}[x]$, we have $G_1\cong G_0\cong H$, with $u$ corresponding to $s$. Let $J=G_1+v$. Then $J$ is a $G_1$-add since $N_J(v)= N_{G_0}(x)\ne\emptyset$. Moreover, since $v$ is nonadjacent to $u$, we have $N_J(v)\not\in\{\{u\},N_{G_1}[u]\}$. Thus $J$ satisfies the requirement.
\hfill $\Box$\medskip

The next theorem is the prime counterpart of Theorem \ref{thm:addirr} and is stated in a slightly stronger form. Together with Lemma \ref{lem:irradd}, it immediately implies that every claw-free prime irreducible graph $G$ properly containing a nontrivial prime irreducible induced subgraph $H$ with at most one simplicial vertex contains a prime irreducible $H$-add.

\begin{theorem} \label{thm:primeadd-irr}
Let $G$ be a prime, irreducible, claw-free graph. Let $H\prec G$, where $H$ is a nontrivial prime graph with at most one simplicial vertex. Then there exists a prime $H$-add $H^+$ such that $H^+\preceq G$ and $d_{H^+}(v)\ge2$, where $v$ is the added vertex.
\end{theorem}

\noindent{\bf Proof}. We first consider the case in which $H$ has no simplicial vertices. By Lemma \ref{lem:primeadd}, there exists a prime $H$-add $H^+ \preceq G$ with added vertex $v$. If $d_{H^+}(v)=1$, then $v$ extends to a claw in $H^+$, contradicting the assumption that $G$ is claw-free. Hence $d_{H^+}(v)\ge2$, so $H^+$ is the desired graph. For the remainder of the proof, assume that $H$ has a unique simplicial vertex. We prove the theorem by identifying $H_0\subset J\subseteq G$ such that $H_0\cong H$, $J=H_0+v$ is prime, and $d_J(v)\ge2$. 

By Theorem \ref{thm:addirr}, there exist $G_0\subset G_0^+\subseteq G$ such that $G_0\cong H$, $G_0^+$ is a $G_0$-add with added vertex $x$, and $N_{G_0^+}(x)\notin \{\{s\}, N_{G_0}[s]\}$, where $s$ is the unique simplicial vertex of $G_0$. Moreover, $d_{G_0^+}(x)\ge2$, for otherwise the unique member of $N_{G_0^+}(x)$ would not be simplicial in $G_0$, so $x$ could extend to a claw in $G_0^+$, a contradiction. If $G_0^+$ is prime, then $J=G_0^+$ is the desired graph. Hence we assume that $G_0^+$ is not prime. Since $G_0\cong H$ is prime, by Observation \ref{obs:basic}(2), $x$ has a twin $y$ in $G_0^+$. The choice of $x$ implies that $y\ne s$, and thus neither $x$ nor $y$ is simplicial in $G_0^+$. Since $G$ is prime, it has a vertex $z$ adjacent to exactly one of $x$ and $y$, say $x$. As $x$ is not simplicial in $G_0^+$, we must have $|N_G(z)\cap V(G_0^+)|\ge2$. Let $G_0^*=G_0^++z$. Since $d_{G_0^+}(x)\ge2$, the induced path $zxy$ is not a component of $G_0^*$. Thus $(G_0^*, x, y, z)$ is a quadruple and it satisfies the assumptions of Lemma \ref{lem:partition}. In what follows we use the notation and conclusions of that lemma, with $G$ replaced by $G_0^*$. 

Suppose both $a_0$ and $b_{n+1}$ are present. By \ref{lem:partition}(2), $J=G_0^*-a_0$ is a prime $(G_0^*-a_0-b_1)$-add with added vertex $b_1$. Since $b_1$ is adjacent to $a_1$ and $b_{n+1}$, $J$ is the desired graph. It remains to consider the case in which $a_0$ or $b_{n+1}$ is absent. By symmetry, we assume that $a_0$ is absent.

To apply \ref{lem:partition}(3), we first claim that $d_{G_0^*}(a_1)>1$. Otherwise $n=1$, so $a_1b_1=xz$, contradicting either $d_{G_0^*}(x)\ge3$ or $d_{G_0^*}(z)\ge2$. Now by \ref{lem:partition}(3), $G_0^*-b_1$ is a prime $(G_0^*-a_1-b_1)$-add with added vertex $a_1$. If $d_{G_0^*}(a_1)\ge3$ then $J=G_0^*-b_1$ is the desired graph. Hence we assume $d_{G_0^*}(a_1)=2$, which implies $n\le2$. We also observe that $B_1$ is a clique, since otherwise two nonadjacent vertices of $B_1$, together with $a_1$ and $b_1$, would induce a claw.

Suppose that $n=2$. Since $a_1$ is complete to $A_1\cup C$, we must have $A_1=C=\emptyset$. Furthermore, $b_3$ cannot be present, because otherwise both $a_2$ and $b_3$ would be simplicial in $G^*_0-a_1-b_1\cong H$, a contradiction. Then by \ref{lem:partition}(3), $J=G_0^*-a_2$ is a prime $(G_0^*-a_2-b_2)$-add with added vertex $b_2$. Moreover, $d_J(b_2)\ge2$, for otherwise $B_1=\emptyset$, implying that $a_1$ and $b_1$ are twins in $G^*_0-a_2-b_2\cong H$, again a contradiction. Thus $J$ is the desired graph.

Finally, assume that $n=1$. Since $|A_0\cup B_0|\ge3$, the vertex $b_2$ is present. Since $a_1b_1=xz$, $G_0^*-b_1$ is prime (by \ref{lem:partition}(3)), and $G_0^*-z$ has twins $x$ and $y$, we must have $a_1=z$ and $b_1=x$. In $G_0^*-a_1=G_0^*-z$, the vertex $b_1=x$ is a twin of both $b_2$ and $y$. Thus either $b_2=y$ or $b_2$ and $y$ are twins. If $b_2\ne y$, then Observation \ref{obs:basic}(1) implies that they remain twins in $G_0^*-a_1-b_1\cong H$, a contradiction. Hence, $b_2=y$. Consequently, $C\ne\emptyset$, for otherwise both $x$ and $y$ would be simplicial in $G^*_0-z=G_0^+$, a contradiction. Since $N_{G_0^*}(a_1) =\{b_1\}\cup A_1\cup C$ and $d_{G_0^*}(a_1)=2$, we have $|C|=1$ and $A_1=\emptyset$. Let $c$ be the unique vertex in $C$. Then $c$ is anticomplete to $D$, because if $c$ has a neighbor $d\in D$, then $\{c,d,a_1,b_2\}$ would induce a claw. Furthermore, $c$ has a neighbor in $B_1$. Otherwise, since $A_1=\emptyset$ and $c$ is anticomplete to $B_1\cup D$, moving $c$ from $C$ to $A_0$ and relabeling it as $a_2$ yields a partition satisfying (P1) and (P2) with larger $|A_0\cup B_0|$, a contradiction. Hence, $B_1\ne\emptyset$. Now $J=G_0^*-b_2$ is a $(G_0^*-a_1-b_2)$-add with added vertex $a_1$. Since $N_J(a_1)=\{c,b_1\}$, we have $d_J(a_1)=2$. Moreover, $c$ and $b_1$ both have degree exceeding 2 in $J$, so $a_1$ is a twin of neither. By Observation \ref{obs:basic}(2), $J$ is prime. Thus $J$ is the desired graph.
\hfill$\Box$\medskip

In both theorems of this subsection, the assumption that $H$ has at most one simplicial vertex is essential for guaranteeing the existence of the desired $H$-add. However, an $H$-add may introduce additional simplicial vertices, so the construction cannot be iterated. Consequently, unlike the splitter theorems in the previous subsections, these results do not yield a chain of (prime) irreducible graphs from $H$ to $G$. To address this issue, we introduce a modified add operation that can be iterated and thus leads to a desired chain. 

For any graph $H$, we define an $H$-add$_2$ as follows. If $H$ has at most one simplicial vertex, then an $H$-add$_2$ is an $H$-add in which the added vertex has degree at least two; otherwise, an $H$-add$_2$ is simply an $H$-add. Equivalently, every $H$-add$_2$ is an $H$-add, and the only $H$-adds that are not $H$-add$_2$s are those in which $H$ has at most one simplicial vertex and the added vertex has degree~1.

The following theorem is an immediate consequence of Lemma \ref{lem:primeadd} and Theorem~\ref{thm:primeadd-irr}. It implies that every prime irreducible claw-free graph can be constructed from any of its nontrivial prime connected induced subgraphs by a sequence of add$_2$ operations. Although the intermediate graphs need not all be irreducible, they remain prime and connected. Moreover, each intermediate graph is irreducible whenever its predecessor is irreducible and has at most one simplicial vertex.

\begin{theorem} \label{thm:chain-irr}
Let $G$ be prime, irreducible, and claw-free. Let $H\prec G$, where $H$ is nontrivial and prime. Then there exists a sequence $H \cong G_0\prec G_1\prec \cdots \prec G_k\cong G$ of prime graphs such that each $G_i$ $(i=1,...,k)$ is a $G_{i-1}$-add$_2$.
\end{theorem}

The add operation can be modified similarly to obtain a chain version of Theorem~\ref{thm:addirr}. We omit the details, as the argument is entirely analogous.

\section{Finite test theorems}

Let $\cal F$ be a family of forbidden graphs and let $H$ be a prime, connected, $\cal F$-free, nontrivial graph. We are interested in characterizing all prime, connected, $\cal F$-free graphs that contain an induced $H$. By Theorem \ref{thm:chain-connp}, every such graph can be constructed from $H$ through a sequence of add operations. To formalize this process, let $\mathcal H_0 =\{H\}$ and, for each $i\ge1$, let $\mathcal H_i=\{G: G$ is a prime, $\cal F$-free $G'$-add for some $G'\in\mathcal H_{i-1}\}$. Then $\bigcup_{i\ge0} \mathcal H_i$ is the set of graphs we are looking for. Since each $\mathcal H_i$ is finite, it can be explicitly generated. In practice, the first few iterations often reveal a structural pattern. Suppose we identify a candidate class $\cal H$ such that, for some $k$, the graphs in $\cal H$ of order at most $|H|+k$ are exactly graphs in $\mathcal H_0\cup \cdots \cup \mathcal H_k$. This raises a natural question: is there a threshold $k$ for which this local agreement implies the global equality $\mathcal H = \bigcup_{i\ge0}\mathcal H_i$? In general, such a threshold need not exist. However, in the following, we identify several scenarios, under different connectivity settings, where such a threshold does exist. 

In this section, we frame and analyze problems in terms of complement graphs for convenience. Accordingly, if a concept X applies to $\ov G$, then we refer to its counterpart in $G$ as {\it coX}. For instance, two vertices are {\it cotwins} in $G$ if they are twins in $\ov G$, $G$ is {\it coprime} if $\ov G$ is prime, and the complement of a $\ov G$-add is a $G$-{\it coadd}. We will use similar terms without further definition.

We also use the following terminology throughout this section. Let $\cal F$ denote a finite set of coprime graphs and let $f = \max\{|F| : F \in \cal F\}$. We define $f=0$ if $\cal F=\emptyset$. For any class $\cal S$ of graphs and any two integers $m\le n$, let $\mathcal S_{m,n}$ denote the set of graphs $G\in \cal S$ with $m\le |G|\le n$. We abbreviate $\mathcal S_{n,n}$ as $\mathcal S_n$. For any set $\cal X$ of graphs and any $i\in \{0,1,2\}$, let 

\makebox[11mm]{} cpadd$_i(\mathcal X, \mathcal F) = \{G : G$ is an $\cal F$-free coprime $X$-coadd$_i$ of some $X \in \cal X\}$,

\noindent where coadd$_0$ is coadd. Iteratively, let cpadd$_i^0(\cal X,F) = X$ and for all integers $k\ge1$, let 

\makebox[11mm]{} cpadd$_i^k(\cal X, F) =$ cpadd$_i$(cpadd$_i^{k-1}(\cal X, F), F)$.

To state our first two main results of this section we need some more definitions. Let $R$ be a fixed and possibly empty set of vertices.  A graph $G$ is called \emph{$R$-rooted} if $R\subset V(G)$. Two $R$-rooted graphs are {\it isomorphic} (or $R$-{\it isomorphic}) if they admit a graph isomorphism that is the identity on $R$. Let $\mathcal A$ (which could be empty) and $\mathcal B$ be two finite sets of $R$-rooted graphs where $G-R$ is connected for every $G\in\cal A\cup B$. If $\mathcal A\neq\emptyset$, let $a=\max\{|A-R|:A\in\mathcal A\}\ge1$, and if $\mathcal A=\emptyset$, let $a=0$. Let $b=\max\{|B-R|:B\in\mathcal B\}\ge2$. 

We define $\sigma(\mathcal A,\mathcal B)$ to be the set of $R$-rooted graphs $G$ for which the components of $G-R$ can be ordered as $G_0,G_1,...,G_t$, where $t\ge 1$ and $G_0$ may be $\oslash$ (indicating the absence of the term $G_0$), such that $G_i+R$ is $R$-isomorphic to a graph in $\cal B$ for every $i\in [t]$, and, if $G_0\ne\oslash$, then $G_0+R$ is $R$-isomorphic to a graph in $\cal A$. In addition, we also require that, if $|G_i|=|G_j|=1$ then $G_i+R$ and $G_j+R$ are not $R$-isomorphic. We write $G = \sigma(G_0, G_1, ..., G_t)$ whenever such an ordering is fixed. The resulting decomposition need not be unique since, when $\cal A\cap B\ne\emptyset$, the role of $G_0$ may be played by another $G_i$. If $R=\emptyset$, note that $\sigma(\cal A,B)$ is the family of graphs that are disjoint unions of at least one graph from $\cal B$ (with repetitions allowed), together with at most one graph from $\cal A$, and that have at most one isolated vertex.
Let $\mathcal S = \sigma(\cal A, B)$, $d=\max\{a-1,b\}$, $r_0=\max\{1,|R|\}$, and $r_2=\max\{2,|R|\}$. The following are the first two main results of this section.

\begin{theorem}\label{thm:cpadd}
Let integers $n$ and $p$ satisfy $p\ge b+1$, $n\ge |R|+a+r_0b+d$, and $n+p\ge |R|+a+$ $b\max\{2,f\}$. Suppose that $\ov G$ is connected for every $G\in\mathcal S_{n+1}$. If $\mathrm{cpadd}^k(\mathcal S_n,\mathcal F)=\mathcal S_{n+k}$ holds for every $k\in [p]$, then it holds for every $k>p$.
\end{theorem}

\begin{theorem}\label{thm:cpadd2}
Let integers $n$ and $p$ satisfy $p\ge b+1$, $n\ge |R|+a+r_2b+d$, and $n+p\ge |R|+a+bf$. Suppose that $\overline{K_{1,3}}\in\cal F$ and $\ov G$ is irreducible for every $G\in\mathcal S_{n+1,n+p}$. If $\mathrm{cpadd}_2^k(\mathcal S_n,\mathcal F)=\mathcal S_{n+k}$ holds for all $k\in [p]$, then it holds for all $k>p$.
\end{theorem}

The third main result of this section is Theorem \ref{thm:ftt2}, which addresses double connectivity for a different class of graphs. We defer this result until after the proofs of Theorems \ref{thm:cpadd} and \ref{thm:cpadd2}. Because these two proofs overlap significantly, we isolate their common arguments into several lemmas, which we then use to prove the two theorems. Let $G=\sigma(G_0,G_1,\ldots, G_t)\in\cal S$. If $H=\sigma(G_0,G_{i_1},...,G_{i_s})$ for some nonempty subset $\{i_1, ..., i_s\}\subseteq [t]$, then we write $H\sqsubseteq G$. Note that $H\sqsubseteq G\in\cal S$ implies $H\in\cal S$.

\begin{lemma}\label{X<H}
Suppose $c\ge 1$, $p\ge b$, and $n+p\ge |R| + a + bc$. Let $G\in \cal S$ with $|G|\ge n+1$, and let $X\subseteq V(G)$ with $|X|\le c$. Then there exists $H\sqsubseteq G$ satisfying $X\subseteq V(H)$ and $H\in\mathcal S_{n+1,n+p}$.
\end{lemma}

\noindent {\bf Proof}. Reorder $G_1, ..., G_t$ so that those containing a vertex of $X$ come first. Let $s\ge 1$ be the smallest index such that $H_s=\sigma(G_0,G_1,...,G_s)$ contains $X$. Then $|H_s|\le |R|+a + bc \le n+p$. If $|H_s|\ge n+1$, then $H_s$ satisfies the requirements. If $|H_s|\le n$, since $|G|\ge n+1$, there exists a smallest index $r$ such that $H_r=\sigma(G_0,G_1,...,G_r)$ has order at least $n+1$. Since $r\ge s$ and $r$ is minimal, $X\subseteq V(H_s)\subseteq V(H_r)$ and $|H_r|\le n+b\le n+p$. Hence $H_r$ satisfies the requirements.\hfill$\Box$ 

\begin{lemma}\label{lem:S1}
Suppose $p\ge b$, $n+p\ge |R|+a+b\max\{2,f\}$, and every graph in $\mathcal S_{n+1,n+p}$ is coprime and $\cal F$-free. Then every $G\in\cal S$ with $|G|\ge n+p+1$ is coprime and $\cal F$-free. In addition, there exists $H\in \mathcal S_{n+1,n+p}$ satisfying $H\sqsubseteq G$.
\end{lemma}

\noindent {\bf Proof}. Let $G\in\mathcal S_{n+k}$, where $k\ge p+1$. We first prove that $G$ is coprime. Suppose, to the contrary, that $G$ contains cotwins $x$ and $y$. Applying Lemma \ref{X<H} with $X=\{x,y\}$ and $c=2$, we obtain $H\sqsubseteq G$ such that $X\subseteq V(H)$ and $H\in\mathcal S_{n+1,n+p}$. By Observation \ref{obs:basic}(1), $x$ and $y$ remain cotwins in $H$, contradicting the hypothesis that every graph in $\mathcal S_{n+1,n+p}$ is coprime.

Next we prove that $G$ is $\cal F$-free. Suppose, to the contrary, that $G$ contains an induced $F\in\cal F$. Applying Lemma \ref{X<H} with $X=V(F)$ and $c=f$, we obtain $H\sqsubseteq G$ with $F\subseteq H\in\mathcal S_{n+1,n+p}$, contradicting the hypothesis that every graph in $\mathcal S_{n+1,n+p}$ is $\cal F$-free.

Finally, applying Lemma \ref{X<H} with $X=\emptyset$ and $c=1$, we obtain $H\sqsubseteq G$ such that $H\in\mathcal S_{n+1,n+p}$, which satisfies the last requirement of the lemma. \hfill$\Box$

\begin{lemma}\label{lem:S2}
Let $i\in\{0,2\}$, $p\ge b+1$, $n\ge |R|+a+r_ib+d$, and $\overline{K_{1,3}}\in\cal F$ if $i=2$. Suppose that $\mathrm{cpadd}_i^k(\mathcal S_n,\mathcal F)=\mathcal S_{n+k}$ for every $k\in [p]$. If $G\in \cal S$ with $|G|\ge n+1$, and $G^+$ is a coprime $\cal F$-free $G$-coadd$_i$, then $G^+\in\cal S$.
\end{lemma}

\noindent{\bf Proof}. Since $G\in \cal S$, we may write $G=\sigma(G_0, G_1, ..., G_t)$. Let $v$ be the added vertex in $G^+$ and let $X=\{x\in R: x$ and $v$ are cotwins in $G^+[R\cup V(G_0)\cup \{v\}]\}$. Since $G^+$ is coprime, for each $x\in X$, we may choose a vertex $x'$ adjacent to exactly one of $x$ and $v$. Let $X'=\{x': x\in X\}$. Since $G^+$ is a $G$-coadd$_i$, the set $M=V(G)-N_{G^+}(v)$ of non-neighbors of $v$ is nonempty. We define $M_i$ as follows. Let $M_0$ be a one-element subset of $M$, giving preference to vertices of $R$. When $i=2$, we consider two cases. If $|M|\ge 2$, let $M_2$ be a two-element subset of $M$, again giving preference to vertices of $R$. If $|M|=1$, then the definition of coadd$_2$ implies that $\ov G$ has at least two simplicial vertices. In addition, since $\ov G \subseteq \ov{G^+}$ is claw-free, the unique member $u$ of $M$ must be simplicial. In this case, let $M_2$ consist of $u$ and another simplicial vertex of $\ov G$. 

Let $\{j_1,...,j_s\}\subseteq [t]$ (where $s=s_i$) be a minimal nonempty set such that $\sigma(G_0,G_{j_1},...,G_{j_s})$ contains $X'\cup M_i$. Without loss of generality, let $j_k = k$ for all $k\in [s]$. Then reorder $G_{s+1}, ..., G_t$ so that those that contain a neighbor of $v$ come first. 

We claim that $s_i\le r_i$. Since the proofs for $i=0$ and $i=2$ follow the same logic, with the former being straightforward, we present only the more involved case $i=2$. 
%By the definition of $M_2$, it is immediate that $\min\{2, |M\cap R|\}\le |M_2\cap R|$. This inequality will be used only in the last two cases below.
In the last two cases below, we will use $\min\{2, |M\cap R|\}\le |M_2\cap R|$, which follows immediately from the definition of $M_2$. 
If $|M_2\cap R|\ge 2$, then $M_2\subseteq R$, and hence $s_2\le \max\{1,|X'|\}\le \max\{1, |X|\}\le \max\{1,|R|\}$. If $|M_2\cap R|=0$, then $|M\cap R|=0$, implying $X = X' = \emptyset$, since $X\subseteq M\cap R$. Thus $s_2\le |M_2|=2$. Finally, if $|M_2\cap R|=1$, then $|X|\le |M\cap R|\le 1$, so $|X'|\le 1$, and therefore $s_2 \le |X'|+ |M_2-R|\le 2$. In conclusion, $s_2\le \max\{2,|R|\}=r_2$, proving the claim.

Since $|G_0|\le a$ and $|G_j|\le b$ for all $j\in [s]$, we have $|\sigma(G_0,G_1,...,G_s)|\le |R|+a+ r_ib\le n-b$. As $|G|\ge n+1$, there exists a minimal index $h=h_i>s_i$ such that $H_i = \sigma(G_0, G_1, ..., G_h)$ satisfies $|H_i|\ge n+1$. By the minimality of $h$, $|H_i|\le n+b\le n-1+p$. Hence, $H_i\sqsubseteq G$ and $H_i\in \mathcal S_{n+1, n+p-1}$. Moreover, since $M_i\subseteq V(H_i)$, it is routine to verify that $H_i^+=H_i+v$ is an $H_i$-coadd$_i$.

Next, we show that $H_i^+$ is coprime. Suppose otherwise. Since $H_i$ is coprime, Observation \ref{obs:basic}(2) implies that $v$ has a cotwin $v'$ in $H_i^+$. As $X'\subseteq V(H_i)$, Observation \ref{obs:basic}(1) gives $v'\notin R$. Thus $v'\in V(G_j)$ for some $j\in\{0,1,...,h\}$, forcing every neighbor of $v$ in $H_i^+$ to lie in $R\cup V(G_j)$. Hence $j\le s+1$, for otherwise the ordering of $G_{s+1}, ..., G_t$ would be violated. Moreover, $h\ge s+2$ because $|H_i| -|\sigma(G_0,G_1,...,G_s)|\ge (n+1) - (n-b) =b+1$. Therefore, $v$ has no neighbors in $G_h$, and hence none in $G_h\cup \cdots  \cup G_t$. Thus $v$ and $v'$ are also cotwins in $G^+$, contradicting the choice of $G^+$.

Since $H_i^+$ is a coprime $\cal F$-free (as $H_i^+\subseteq G^+$) $H_i$-coadd$_i$, we have $H_i^+\in \text{cpadd}_i(\mathcal S_{n+1,n+p-1},\mathcal F) =\mathcal S_{n+2,n+p}$. Let $H_i^+ = \sigma(J_0,J_1,...,J_k)$ ($k=k_i$) and $j$ be the index with $v\in V(J_j)$. We claim that $v$ is anticomplete to $G_h$. If not, by the ordering of $G_{s+1}, ..., G_t$, $v$ has a neighbor in each of $G_{s+1},...,G_h$. Hence, $d+1\ge |J_j|\ge |H_i^+|-|\sigma(G_0,G_1,...,G_s)|\ge 1+(n+1)-(|R|+a+r_ib)$, contradicting the hypothesis on $n$. This proves the claim. Therefore, $v$ is anticomplete to $G_h\cup \cdots  \cup G_t$, and $G^+=\sigma(J_0,J_1,...,J_k, G_{h+1}, ..., G_t)$, implying $G^+\in \cal S$ and completing the proof of the lemma.\hfill$\Box$\medskip

\noindent{\bf Proof of Theorem \ref{thm:cpadd}}. With the preceding preparations in place, the proof of this theorem is straightforward. Let $k>p$. The inclusion $\mathrm{cpadd}^k(\mathcal S_n,\mathcal F)\subseteq \mathcal S_{n+k}$ follows immediately from Lemma \ref{lem:S2}. For the reverse inclusion, let $G\in\mathcal S_{n+k}$. By Lemma \ref{lem:S1}, $G$ is coprime and $\cal F$-free, and there exists $H\sqsubseteq G$ with $H\in\mathcal S_{n+j}$ for some $j\in [p]$. Since $\mathcal S_{n+j}=\text{cpadd}^{j-1}(\mathcal S_{n+1},\cal F)$, $H$ is obtained from some $X\in \mathcal S_{n+1}$ via a sequence of coadd operations. Because $\ov X$ is connected, this implies that $\ov H$ is connected. Moreover, $H$ is coprime as $H\in \text{cpadd}^j(\mathcal S_n,\cal F)$, and $\ov G$ is connected as $\ov G-\ov H$ is complete to $\ov H-R$ in $\ov G$. Therefore, by Theorem \ref{thm:chain-connp}, $G$ can be obtained from $H$ via a sequence of coadd operations. Thus $G\in\text{cpadd}^{k-j}(\mathcal S_{n+j},\mathcal F) = \text{cpadd}^k(\mathcal S_n,\cal F)$, completing the proof.\hfill$\Box$\medskip

\noindent{\bf Proof of Theorem \ref{thm:cpadd2}}. The proof is very similar to that of the last theorem. Let $k > p$. As before, the inclusion $\mathrm{cpadd}_2^k(\mathcal S_n,\mathcal F)\subseteq \mathcal S_{n+k}$ follows immediately from Lemma \ref{lem:S2}. For the reverse inclusion, let $G \in \mathcal S_{n+k}$. By Lemma \ref{lem:S1}, $G$ is coprime and $\cal F$-free, and there exists $H\sqsubseteq G$ with $H\in \mathcal S_{n+j}$ for some $j\in [p]$. 

We first prove that $\ov G$ is irreducible. Note that $\ov H$ is irreducible since $H\in \mathcal S_{n+j}$. In particular, by Observation \ref{obs:basic}(3), $\ov H$ is 2-connected. Since $|H-R|= n+j - |R|\ge r_2b+d\ge 2\cdot 2+2=6$ and $\ov G -\ov H$ is complete to $\ov H-R$ in $\ov G$, it follows that $\ov G$ is 2-connected, and thus neither disconnected nor a simplicial sum. It remains to show that $\ov G$ is not a simplicial join. Suppose to the contrary that $\ov G$ is the simplicial join of $G_1$ and $G_2$ over their respective simplicial vertices $v_1$ and $v_2$. Let $J=G_1-N_{G_1}[v_1]$. By Observation \ref{obs:basic}(4), we may assume $\ov H \cap J = \emptyset$. Note that $J \ne \emptyset$; otherwise, the vertices of $N_{G_1}(v_1)$ would be twins of $\ov G$, contradicting that $G$ is coprime. Since $J\subseteq \ov G-\ov H$, $J$ is complete to $\ov H-R$, which implies $V(\ov H -R) \subseteq N_{G_1}(v_1)$, and thus, by Observation \ref{obs:basic}(1), the vertices of $\ov H-R$ are twins of $\ov H$, contradicting that $H$ is coprime. Therefore, $\ov G$ is irreducible.

In summary, $\ov G$ is prime, irreducible, and claw-free, while $\ov H\subseteq \ov G$ is prime, irreducible, and nontrivial. By Theorem \ref{thm:chain-irr}, $G$ can be constructed from $H$ by a sequence of coadd$_2$ operations. Therefore, $G\in\text{cpadd}^{k-j}_2(\mathcal S_{n+j},\mathcal F) = \text{cpadd}^k_2(\mathcal S_n,\cal F)$, which proves the theorem.\hfill$\Box$\medskip

Our third main result of this section involves bipartite graphs. Since such graphs could admit more than one partition into two stable sets, which would introduce ambiguity in our discussions, we assume that each bipartite graph comes with a $\{0,1\}$-labeling (a 2-coloring) of its vertices such that vertices with the same label form a stable set. We refer to these as {\it labeled graphs} and denote them by $G = (X, Y, E)$, where $X$ and $Y$ are the two color classes containing the vertices labeled $0$ and $1$, respectively. The disjoint union of two labeled graphs is defined naturally. Additionally, we define the {\it bipartite complement} of a labeled graph $G = (X, Y, E)$ as the labeled graph $b(G) = (X, Y, F)$, where, for every $x\in X$ and $y\in Y$, $xy\in F$ if and only if $xy\notin E$. Note that $b(b(G))=G$.

For our third main theorem, we redefine $\cal B$ and $\cal S$ as follows. Let $\cal B$ be a finite set of labeled connected coprime graphs such that every connected coprime induced subgraph of a graph in $\cal B$ also belongs to $\cal B$. We assume $b = \max\{|B| : B \in \mathcal{B}\} > 1$. Let $\tau(\cal B)$ be the set of labeled graphs  $G$ such that $b(G)$ is a disjoint union of copies of graphs in $\cal B$, with the additional condition that in $b(G)$, no two isolated vertices have the same label. Note that if $G\in\tau(\cal B)$ and $b(G)$ is the disjoint union of $B_i=(X_i,Y_i,E_i)$ for  $i\in [k]$, then $\ov G$ is the graph obtained from $B_1\cup \cdots  \cup B_k$ by turning both $X_1\cup \cdots \cup X_k$ and $Y_1\cup \cdots  \cup Y_k$ into cliques. In particular, $\ov G$ is claw-free. Let $\mathcal S = \tau(\cal B)$. We first establish a few properties of $\cal S$. 

\begin{lemma}\label{lem:b(G)}
Let $G=(X,Y,E)$ be a labeled graph. Then the following hold. 
\begin{enumerate}[label=(\arabic*)]
\parskip -4pt
\item Any two vertices of the same label are cotwins in $G$ if and only if they are cotwins in $b(G)$. 
\item If $G\in \cal S$ and $b(G)$ is disconnected, then $G$ is coprime. 
\item If $H\subseteq G\in \cal S$ and $H$ is coprime, then $H\in\cal S$.
\end{enumerate} 

\end{lemma}

\noindent{\bf Proof}. We omit the proof of (1) since it is straightforward. 

To prove (2), suppose for a contradiction that $G$ has cotwins $u$ and $v$. Without loss of generality, assume that $u\in X$. If $v\in Y$, then the conditions $N_G(u) \subseteq Y$, $N_G(v) \subseteq X$, and $N_G(u)=N_G(v)$ imply $N_G(u) = N_G(v) = \emptyset$. Hence, in $b(G)$, $u$ is complete to $Y$ and $v$ is complete to $X$, so $b(G)$ is connected, contradicting the hypothesis of (2). Therefore, $v\in X$. By (1), $u$ and $v$ are cotwins in $b(G)$. Since each component of $b(G)$ belongs to $\cal B$ and each graph in $\cal B$ is coprime, $u$ and $v$ must lie in distinct components of $b(G)$. Thus, as cotwins in $b(G)$ with the same label, $u$ and $v$ must both be isolated vertices. This contradicts the definition of $\cal S$, and hence completes the proof of (2).

To prove (3), note that since $H$ is coprime, no two isolated vertices in $b(H)$ may have the same label. By the definition of $\cal S$, it suffices to show that every component $H_0$ of $b(H)$ belongs to $\cal B$. Because $H\subseteq G$, we have $b(H) \subseteq b(G)$, implying $H_0\subseteq G_0$ for a component $G_0$ of $b(G)$. Since $G \in \cal S$, we have $G_0\in \cal B$. By the definition of $\cal B$, to prove $H_0\in \cal B$, we need only show that $H_0$ is coprime. Suppose for contradiction that $H_0$ contains cotwins $u$ and $v$. If $u$ and $v$ have different labels, since $H_0$ is connected and bipartite, $u$ has a neighbor not adjacent to $v$, a contradiction. If $u$ and $v$ have the same label, they are cotwins in the entire $b(H)$, which by (1) implies they are cotwins in $H$, contradicting that $H$ is coprime. Thus, $H_0$ is coprime, completing the proof.\hfill$\Box$\medskip

The following is the third main result of this section which is about doubly connected graphs.

\begin{theorem} \label{thm:ftt2}
Let integers $n$ and $p$ satisfy $n\ge 2b+1$, $p\ge b$, and $n+p\ge f$. If $cpadd_1^k(\mathcal S_n, \mathcal F) = \mathcal S_{n+k}$ holds for all $k\in [p]$, then the same holds for all $k>p$.
\end{theorem}

\noindent{\bf Proof}. This proof follows a similar argument to that of theorems \ref{thm:cpadd} and \ref{thm:cpadd2}.  
We first prove that all graphs in $\cal S$ are $\cal F$-free. Suppose otherwise; then some $F\in \cal F$ is contained in a graph of $\cal S$. Since $F$ is coprime, Lemma \ref{lem:b(G)}(3) implies that $F\in \cal S$. As $|F|\le f\le n+p$ and all graphs in $\mathcal S_{n+1,n+p}$ are $\cal F$-free, we have $|F|\le n$. Because $K_2\in\cal B$, by the definition of $\cal S$, the disjoint union of $b(F)$ with any number of $K_2$ is equal to $b(G)$ for some $G\in\cal S$. If we choose $G$ with $n< |G|\le n+2$, then $F\subseteq G\in\mathcal S_{n+1,n+p}$, yielding a contradiction.

Next, we show that every graph $G=(X,Y,E)\in\cal S$ with $|G|\ge n$ is doubly connected. By the definition of $\cal S$, $b(G)$ has at most two isolated vertices. Since $|G|\ge n\ge3$, $b(G)$ contains a nontrivial component $B_1=(X_1,Y_1,E_1)$. As $\ov G$ is obtained from $b(G)$ by turning $X$ and $Y$ into cliques, the edges of $B_1$ connect the two cliques, and hence $\ov G$ is connected. 
It remains to show that $G$ is connected. Let $X_2=X-X_1$ and $Y_2=Y-Y_1$. Since $|G|\ge n\ge b+2$ and no two isolated vertices of $b(G)$ have the same label, both $X_2$ and $Y_2$ are nonempty. Moreover, $G[X_1\cup Y_2]$ and $G[X_2\cup Y_1]$ are complete bipartite graphs. Thus, it suffices to show that $G[X_1\cup Y_1]\cup G[X_2\cup Y_2]$ contains an edge. Suppose otherwise. Then $b(G[X_1\cup Y_1])$ and $b(G[X_2\cup Y_2])$ would be complete bipartite graphs, and hence the only two components of $b(G)$. Since all components of $b(G)$ are in $\cal B$ and thus are coprime, we obtain $|X_1|=|Y_1|=|X_2|=|Y_2|=1$, which implies $|G|=4$, contradicting the hypothesis that $|G|\ge n\ge 2b+1\ge 5$. Therefore, $G$ is connected and hence $G$ is doubly connected.

We have shown that all graphs in $\mathcal S_{n+k}$ with $k\ge0$ are $\cal F$-free and doubly connected. To apply Theorem \ref{thm:primedc}, we need a few more observations. For any $G\in \cal S$ with $|G|\ge n$, Lemma \ref{lem:b(G)}(2) implies that $G$ is coprime since $n>b$. Now, let $G\in\mathcal S_{n+k}$ with $k>p$. By successively deleting from $G$ the vertex sets of components of $b(G)$, the first resulting graph of order $\le n+p$ belongs to $\mathcal S_{n+1,n+p}$ because $p\ge b$. Hence, $G$ contains a graph in $\mathcal S_{n+1,n+p}$. We further observe that $\ov G\notin\cal T$ (see Figure \ref{fig:Tn}). Indeed, if $\ov G\in\cal T$, then $b(G)$ would contain a component of order $|G|-2$, which would imply $b\ge |G|-2 > n-2 \ge 2b-1 > b$, a contradiction. Since $\ov G$ is claw-free, Theorem \ref{thm:primedc} ensures that $G$ can be constructed from a graph in $\mathcal S_{n+1,n+p}$ via a sequence of coadd$_1$ operations. Therefore, cpadd$_1^k(\mathcal S_n, \mathcal F) \supseteq \mathcal S_{n+k}$ holds for all $k\ge0$.

It now remains to prove the reverse inclusion, which will occupy the remainder of this proof. Suppose, for a contradiction, that $k\ge p$ is the least integer such that cpadd$_1^{k+1}(\mathcal S_n, \mathcal F) \nsubseteq \mathcal S_{n+k+1}$. Then there exist a graph $G\in \mathcal{S}_{n+k}$ and a coprime $\cal F$-free $G$-coadd$_1$ $G^+$, with added vertex $v$, such that $G^+ \notin \mathcal{S}_{n+k+1}$. Let $N=N_{G^+}(v)$ and $M=V(G)-N$. By the definition of add$_1$, $M\ne\emptyset\ne N$.

We first prove that if $G^+$ is not bipartite, then there exists $G_0=(X_0,Y_0,E_0)\subseteq G$ such that $|G_0|< n$, every component of $b(G_0)$ is a component of $b(G)$, and $N\cap X_0$, $N\cap Y_0$, and $M\cap V(G_0)$ are all nonempty. If some component $B_1$ of $b(G)$ has both color classes intersecting $N$, then, for any component $B_2$ meeting $M$, $G_0 = b(B_1 \cup B_2)$ is a desired graph as $|G_0|\le 2b< n$. Thus, we may assume that in every component of $b(G)$, at most one color class meets $N$. Let $B_i=(X_i,Y_i,E_i)$ for $i=1,2$ be components of $b(G)$ where $X_1$ and $Y_2$ intersect $N$. If $Y_1\cup X_2\ne\emptyset$, then $G_0=b(B_1\cup B_2)$ is a desired graph. If $Y_1\cup X_2=\emptyset$, then $|X_1|=|Y_2|=1$ as $B_1,B_2$ are connected. Now for any component $B_3$ that intersects $M$, $G_0=b(B_1\cup B_2\cup B_3)$ is a desired graph, since $|G_0|\le 2+b\le 2b< n$. This proves the existence of $G_0$. 

We now prove that $G^+$ is bipartite. Suppose otherwise. Since $|G|=n+k\ge n$, $G_0$ can be extended within $G$ to a minimal graph $H$ with $|H|\ge n$ by adding components of $b(G)$ to $b(G_0)$. By the minimality, $|H|\le n-1+b\le n-1+p$, which means $H\in\mathcal S_{n+j}$ for some $j\in\{0,1,...,p-1\}$. Let $H^+= H + v$. The choice of $G_0$ ensures that $H^+$ is a non-bipartite $H$-coadd$_1$. Since $H$ is bipartite, $v$ has no cotwin in $H^+$. Then, by Observation \ref{obs:basic}(2), $H^+$ is coprime. Moreover, $H^+$ is $\cal F$-free as $H^+\subseteq G^+$. Hence, as $j< p$, we have $H^+\in \text{cpadd}_1(\mathcal S_{n+j},\mathcal F)=\mathcal S_{n+j+1}$, contradicting the fact that $H^+$ is non-bipartite. Therefore, $G^+$ must be bipartite.

Since $G^+$ is bipartite, all vertices in $N$ must have the same label, say 0. Then $v$ has label 1. We order the components of $b(G)$ as follows. Let $B_1$ be a component of $b(G)$ that contains a vertex of $N$. We then order the remaining components $B_2$, ..., $B_t$ of $b(G)$ so that those containing a 0-labeled vertex of $M$ appear first. Since $|G| = n+k>n$ and $|B_1|\le b<n$, there is a minimum index $h$ such that $H=b(B_1\cup \cdots  \cup B_h)$ satisfies $|H|\ge n$. The minimality of $h$ implies $|H|< n+b\le n+p$, and hence $H\in\mathcal S_{n+j}$ for some $j\in\{0,1,...,p-1\}$. In addition, since $|H|\ge n\ge 2b+1$, we have $h\ge3$.

Let $H^+= H +v$. Then $H^+$ is an $H$-coadd$_1$ because $B_1$ contains a vertex of $N$ and $B_1\cup B_2$ contains a vertex of $M$; otherwise, $B_1\cup B_2$ would consist of two 0-labeled isolated vertices of $b(G)$, violating $G\in\cal S$. We claim that $H^+$ is coprime. Otherwise, since $H$ is coprime, Observation \ref{obs:basic}(2) implies that $v$ has a cotwin $u$ in $H^+$. As $H$ is connected by the second paragraph of this proof, $u$ and $v$ must have the same label, and thus Lemma \ref{lem:b(G)}(1) implies that they are cotwins in $b(H^+)$. Hence, if $u\in V(B_i)$, then $N_{b(H^+)}(v) \subseteq V(B_i)$. Since $N_{b(H^+)}(v)$ consists precisely of the 0-labeled vertices of $H^+$ that belong to $M$, the ordering of the components of $b(G)$ implies that either both $u$ and $v$ are isolated in $b(H^+)$ or $i\le 2$. In either case, $B_3$ contains no 0-labeled vertices of $M$. By the ordering of the components again, $B_3\cup \cdots  \cup B_t$ contains no 0-labeled vertices of $M$. Therefore, $N_{b(G^+)}(v)=N_{b(H^+)}(v)$. It follows that $u$ and $v$ are cotwins in $b(G^+)$ and hence, by Lemma \ref{lem:b(G)}(1), cotwins in $G^+$, contradicting the assumption that $G^+$ is coprime. This proves the claim. 

Now $H^+$ is a coprime $\cal F$-free $H$-coadd$_1$. It follows that $H^+\in \text{cpadd}_1(\mathcal S_{n+j},\mathcal F)=\mathcal S_{n+j+1}$, as $j<p$. We claim that $v$ is anticomplete to $B_h$ in $b(H^+)$. Suppose otherwise. By the ordering of the components of $b(H)$, $v$ has neighbors in each of $B_2,...,B_h$ in $b(H^+)$. Hence, the component of $b(H^+)$ that contains $v$ has order at least $1 + |B_2|+\cdots+|B_h| = 1 + |H| - |B_1|\ge 1+ n -b>b$. This contradicts $H^+\in\mathcal S_{n+j+1}\subseteq \cal S$, proving the claim. The ordering of the components of $b(G)$ then implies that $v$ is anticomplete to $B_{h+1}\cup  \cdots \cup B_t$ in $b(G^+)$. Consequently, $b(G^+)$ is the disjoint union of $b(H^+)$ and $B_{h+1}, ..., B_t$. Since $H^+\in\cal S$, all components of $b(G^+)$ belong to $\cal B$. In addition, since $G^+$ is coprime, no two isolated vertices of $b(G^+)$ can have the same label. Therefore, $G^+\in\cal S$, contradicting the initial choice of $G^+$. This completes the proof. \hfill$\Box$

\section{Characterizing $\cal F$-free graphs}

\begin{figure}[htbp]
\begin{center}
\includegraphics[scale=0.6]{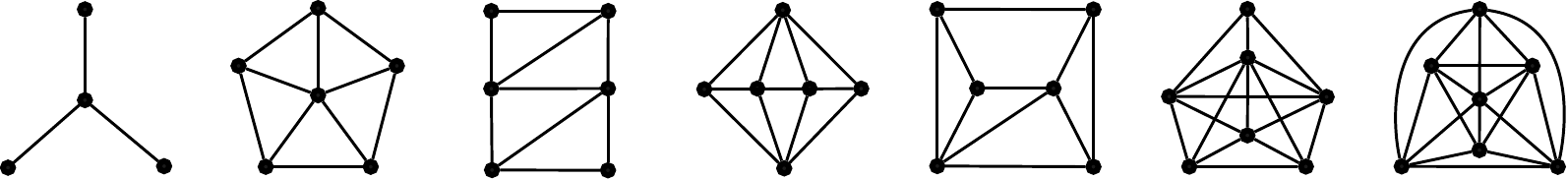}
\caption{Forbidden line graphs $\mathbb L_1,\mathbb L_2,...,\mathbb L_7$.}
\label{fig:line7}
\end{center}
\end{figure}

In this section, we apply the results of the previous sections to the problem of characterizing $\cal F$-free graphs, where $\cal F$ is a set of irreducible prime graphs that includes the claw. We first reduce this problem to the characterization of $\cal F$-free irreducible prime graphs. We then consider a specific case $\mathcal F_0=\{\mathbb L_1,\mathbb L_3\}$, where $\mathbb L_1$ is the claw and $\mathbb L_3$ is another forbidden graph for line graphs (see Figure~\ref{fig:line7}). Using the forbidden-graph characterization of line graphs, we determine all $\mathcal F_0$-free graphs by a computer-assisted finite computation. This example demonstrates how the structural theory developed in this paper can be used to solve such characterization problems. In particular, it provides a framework that separates routine case checking, which can be delegated to a computer, from the essential mathematical proofs.

We begin with a few definitions. {\it Duplicating} a vertex $x$ of a graph $G$ adds a new vertex $y$ that is adjacent to every vertex in $N_G[x]$. It is clear that $x$ and $y$ are twins in the resulting graph. For any family $\mathcal G$ of graphs, let $\mathbf U(\mathcal G)$ denote the class of graphs that are disjoint unions of copies of graphs in $\mathcal G$, let $\mathbf S(\mathcal G)$ denote the class of graphs constructed from copies of graphs in $\mathcal G$ by repeated simplicial sums and simplicial joins, and let $\mathbf D(\mathcal G)$ denote the class of graphs constructed from copies of graphs in $\mathcal G$ by repeatedly duplicating non-simplicial vertices. The following result shows that, in order to characterize all $\cal F$-free graphs, it suffices to characterize the irreducible prime ones. 

\begin{theorem}\label{thm:decomposition}
Let $\cal F$ be a set of irreducible prime graphs that contains $\mathbb L_1$. Then\\ 
\makebox[6mm]{} \{$\cal F$-free graphs\} $= \mathbf U(\mathbf S(\mathbf D(K_1, K_2, K_3, K_4\del e,$ and all irreducible prime $\cal F$-free graphs$)))$.
\end{theorem} 

\noindent{\bf Proof}. 
Here we provide a sketch of the proof. More details can be found in \cite{jake}. 

We proceed in three steps, corresponding to the three operations. First, it is easy to see that, if $F$ is connected and $G$ is the disjoint union of $G_1$ and $G_2$, then $G$ is $F$-free if and only if both $G_1$ and $G_2$ are $F$-free (see 2.1.1 of \cite{jake}). As a result, \{$\cal F$-free graphs\} = $\mathbf U(\cal F$-free connected graphs). 

Our next step is to establish $\{\cal F$-free connected graphs$\} = \{K_1\}\cup \mathbf S(\cal F$-free irreducible graphs). To accomplish this goal, we only need to prove the following.

\begin{enumerate}[label=(\arabic*)]
\parskip -4pt
\item {\it Let $F$ be connected and not a simplicial sum. If $G$ is a simplicial sum of $G_1$ and $G_2$, then $G$ is $F$-free if and only if both $G_1$ and $G_2$ are $F$-free.}
\item {\it Let $F$ be $\mathbb L_1$ or an $\mathbb L_1$-free irreducible graph. If $G$ is a simplicial join of $G_1$ and $G_2$, then $G$ is $F$-free if and only if both $G_1$ and $G_2$ are $F$-free.}
\end{enumerate}

The proof of (1) is straightforward and is therefore omitted (see 5.1.1 of \cite{jake}). Although the proof of (2) follows the same logic (see 5.1.2 of \cite{jake}), we provide a brief sketch because it involves additional case analysis. We only consider the backward implication and we assume, for a contradiction, that $F\subseteq G$. For $i=1,2$, let $v_i$ be the simplicial vertex of $G_i$ over which the join is performed, and set $V_i=V(G_i-v_i)$, $N_i=N_{G_i}(v_i)$, and $U=V(F)$. For each $i$, $U\cap V_i\ne\emptyset$ since otherwise $F\subset G_{3-i}$, and $U\cap N_i\ne\emptyset$ since $F$ is connected. Moreover, $|U\cap V_i|\ge 2$ since otherwise $F\prec G_{3-i}$. Hence $|U \cap N_i| = 1$ for some $i$, say $i=1$, as otherwise $F$ is a simplicial join. Then the unique vertex $u$ of $U\cap N_1$ is a cut vertex of $F$. By Observation \ref{obs:basic}(3), $F\succeq \mathbb L_1$ and thus $F=\mathbb L_1$. In particular, $d_F(u)=3$. Since $|U\cap V_2|\ge2$, $u$ has two neighbors $u_1,u_2\in V_2$. Both $u_1,u_2$ must belong to $N_2$, which is a clique, so $u_1u_2\in E(F)$. Thus $F$ is not a claw, a contradiction. This completes step 2. 

As the final step, it remains to show that \{$\cal F$-free irreducible graphs$\} = \{K_2, K_3\}\cup \mathbf D(K_4\del e$ and $\cal F$-free irreducible prime graphs). To this end (see 5.1.7 of \cite{jake}), we prove the following two statements. Note that the forward inclusion ($\subseteq$) follows immediately from (3). On the other hand, observe that vertex duplication preserves $\cal F$-freeness, as graphs in $\cal F$ are prime. This observation and (4) clearly imply the backward inclusion ($\supseteq$). Therefore, we only need to prove (3) and (4).

\noindent (3) {\it If a claw-free irreducible graph $G$ has twins $x,y$, then the following hold. \\ 
\indent $(a)$ If $d_G(x) \le 2$, then $G\in \{K_2,K_3\}$.\\ 
\indent $(b)$ If $d_G(x) \ge 3$, then $x$ is non-simplicial. Moreover, $G-x$ is irreducible unless $G=K_4\del e$}.\\ 
(4) {\it Duplicating a non-simplicial vertex in a claw-free irreducible graph preserves irreducibility}.

Proving (3a) and the first part of (3b) is straightforward (see 5.1.4 of \cite{jake}): if $d_G(x)=1$, then $G=K_2$ since $G$ is connected; if $d_G(x)=2$, then $G=K_3$ by Observation \ref{obs:basic}(3); if $d_G(x)\ge3$ and $x$ is simplicial, then $G$ is a simplicial join, a contradiction. 
For the second part of (3b) (see 5.1.5 of \cite{jake}), note that $G-x$ is neither disconnected nor a simplicial join, as either would force the same for $G$. Thus, we may assume that $G-x$ is a simplicial sum of $G_1$ and $G_2$ over their common simplicial vertex $v$. Necessarily $v=y$, since otherwise $G$ is a simplicial sum. Moreover, if $d_{G_i}(y)\ge2$ for some $i$, then $G$ is a simplicial join. Hence $d_{G_1}(y)=d_{G_2}(y)=1$, and Observation \ref{obs:basic}(3) gives $G=K_4\backslash e$.

To prove (4), let $G$ be claw-free and irreducible. Let $G^+$ be obtained from $G$ by duplicating a non-simplicial vertex $y$ and let $x$ be the new vertex. Since $G$ is connected, $G^+$ is clearly also connected. Similarly, because $G$ has no cut vertices (by Observation \ref{obs:basic}(3)) and $d_{G^+}(x)\ge2$ (as $d_G(y)\ge1$), $G^+$ has no cut vertices and is therefore not a simplicial sum. Now, assume for a contradiction that $G^+$ is a simplicial join of $G_1,G_2$. For $i=1,2$, let $v_i, V_i, N_i$ be defined as in the proof of (2). Since $G$ is not a simplicial join, $x$ must belong to some $N_i$, say $N_1$, with $|N_1|=2$. Furthermore, because $G$ is not a simplicial sum, $V_1-N_1$ must be empty. Hence $y$ is simplicial in $G$, a contradiction. This completes the proof of (4) and also the theorem.\hfill$\Box$\medskip

As noted earlier, we use computer assistance to perform the tedious but computationally straightforward case analysis. Our implementation is based on SageMath, which provides the necessary graph-theoretic operations. We supplement these capabilities with auxiliary functions and core generation routines. The latter exhaustively generate all candidate graphs and filter them to retain only non-isomorphic $\cal F$-free graphs. A selection of these procedures is described below; the complete list is provided in the Appendix. 

\begin{itemize}
\parskip -4pt
\item {\tt co(G)}: Produces the complement of every graph $\tt g\in G$.
\item {\tt conn(G)}: Outputs the connected members of {\tt G}.
\item {\tt remedy(g, F)}: Generates all $\preceq$-minimal {\tt F}-free {\tt g}-remedies.
\item {\tt Add(G, F, t, p)}: Given $t \in \{0,1,2\}$ and $p \in \{0,1\}$, generates all non-isomorphic {\tt F}-free {\tt g}-add$_t$ graphs for every ${\tt g} \in {\tt G}$. If $p = 1$, only prime graphs are retained.
\item {\tt coAdd(G, F, t, p)} = {\tt co(Add(co(G), co(F), t, p))}.
\end{itemize}

Let $\cal L$ denote the class of line graphs of multigraphs. A classical result of Bermond and Meyer~\cite{line} states that a graph $G$ belongs to $\cal L$ if and only if $G$ contains none of the seven graphs depicted in Figure \ref{fig:line7} as an induced subgraph. However, these seven graphs contribute unequally to carving out the boundary between line graphs and non-line graphs. To illustrate this, let $\mathbb L_2^+$ be the graph obtained from $\mathbb L_2$ by adding a new vertex whose unique nonneighbor in $\mathbb L_2$ is a cubic vertex of $\mathbb L_2$. It is routine to verify (see 2.4.1 of \cite{jake}) that every connected prime $\{\mathbb L_1,\mathbb L_3,\mathbb L_4,\mathbb L_5,\mathbb L_6,\mathbb L_7\}$-free graph must belong to $\mathcal L\cup \{\mathbb L_2,\mathbb L_2^+\}$. Thus, the role of $\mathbb L_2$ is highly localized: it serves only to exclude $\mathbb L_2$ and $\mathbb L_2^+$. This phenomenon is not isolated; in fact, several similar results were obtained in \cite{jake}. In the following, we characterize $\{\mathbb L_1,\mathbb L_3\}$-free graphs, showing that they arise from line graphs, graphs whose complement has girth at least five, and two explicitly described exceptional families.

Let $\cal J$ denote the class of connected coprime graphs of girth at least five. For any positive integers $p$, $q$, and $r$ satisfying $p\ge q\ge r\ge p-1$, let $K_{p,q;r}$ denote the graph obtained from $K_{p,q}$ by deleting a matching of size $r$. Let $\cal K$ consist of $K_{3,3;2}$ and all $K_{p,q;r}$ with $p\ge 4$.

\begin{figure}[htbp]
\begin{center}
\includegraphics[scale=0.6]{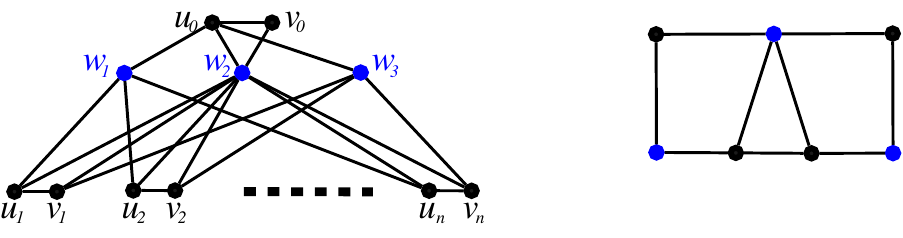}
\caption{Graphs $M_n$ and $M_3^*$.}
\label{fig:Mn}
\end{center}
\end{figure}

For any integer $n\ge3$, let $M_n$ denote the graph (see Figure \ref{fig:Mn}) with vertices $w_1$, $w_2$, $w_3$, $u_0$, $v_0$, $\ldots$, $u_n,v_n$ and edges $u_0w_1,u_0w_2,u_0w_3,u_0v_0,v_0w_2$ and $u_iw_1, u_iw_2, u_iv_i, v_iw_2, v_iw_3$ for all $i\in [n]$. Let $\mathcal M=\{M_n-X: X\subseteq \{u_0,v_0,v_1,u_2\}$ and $n\ge3$\}. Note that the smallest graph in $\cal M$ is $M_3^*=M_3-\{u_0,v_0,v_1,u_2\}$. Let $\mathcal N = \{\ov G: G\in\cal M\}$.

\begin{theorem}\label{thm:13free}
Let $G$ be irreducible and prime. Then $G$ is $\{\mathbb L_1,\mathbb L_3\}$-free if and only if $G\in \cal L\cup N$ or every component of $\ov G$ belongs to $\cal J\cup K$.  
\end{theorem}

\noindent {\bf Proof}. We first characterize $\cal K$ and $\cal N$.

\noindent (5.2.1) {\it Let $H$ be a doubly connected coprime graph that contains an induced $4$-cycle. Then $H$ is $\{K_3, \ov{\mathbb L}_3\}$-free if and only if $H\in \cal K$.} 

To prove the backward implication, we show that every graph $H\in\cal K$ is $\{K_3, \ov{\mathbb L}_3\}$-free, without requiring any additional assumptions on $H$. By the definition of $\cal K$, every induced subgraph of $H$ is obtained from a complete bipartite graph by deleting a matching. Since $K_3$ is non-bipartite and $\ov{\mathbb L}_3=K_{3,3}-E(P_4)$, none can be contained in $H$, which proves the backward implication.

To prove the forward implication, suppose $C=wxyz$ is an induced 4-cycle of $H$. Since $H$ is coprime, it has a vertex $v$ adjacent to exactly one of $w,y$, say $w$. As $H$ is $K_3$-free, $v$ is adjacent to no other vertex of $C$. Likewise, we may choose a vertex $u$ adjacent to $z$ but to none of $w,x,y$. If $uv\notin E(H)$, then $H[u,v,w,x,y,z]\cong \ov{\mathbb L}_3$, a contradiction. Hence $uv\in E(H)$, and consequently $H[u,v,w,x,y,z]\cong K_{3,3;2}$. 

Since $H\succeq K_{3,3;2}$, we have $\ov{H} \succ C_4$. This implies $\ov{H}\notin\cal T$, because every $T_n$ is $C_4$-free (see Figure \ref{fig:Tn}). Let $\mathcal F=\{K_3, \ov{\mathbb L}_3\}$. By Theorem \ref{thm:primedc}, $H\in \bigcup_{k\ge0}\mathbf H_k$, where $\mathbf H_k= {\rm cpadd}_1^k(\{K_{3,3;2}\},\cal F)$. Since $\mathbf H_0=\{K_{3,3;2}\}$ and $\mathbf H_k = {\tt coAdd}(\mathbf H_{k-1}, \mathcal F,1,1)$ for all $k\ge1$, a computer-assisted computation (see Appendix) yields $\mathbf H_1 = \{K_{4,3;3}\}$, $\mathbf H_2 = \{K_{4,4;3},K_{4,4;4}\}$, and $\mathbf H_3 = \{K_{5,4;4}\}$. Therefore, $\mathbf H_k =\mathcal K_{6+k}$ holds for all $k\in \{0,1,2,3\}$. Since, by definition, $\mathcal K=\bigcup_{k\ge0}\mathcal K_{6+k}$, it remains only to show that $\mathbf H_k =\mathcal K_{6+k}$ holds for all $k\ge4$.

To this end, we invoke Theorem \ref{thm:ftt2}. Let $\cal B$ consist of three labeled graphs $(\{x\}, \emptyset,\emptyset)$, $(\emptyset,\{y\},\emptyset)$, and $(\{x\},\{y\},\{xy\})$, and let $\cal S = \tau(B)$. Then $\mathcal S_k= \mathcal K_k$ for $k\ge7$, when the members of $\cal K$ are regarded as labeled graphs (note that each connected bipartite graph can be regarded as a labeled graph in two ways, depending on which color class is labeled 0). Under this choice of $\cal B$ and $\cal F$, we have $b=2$ and $f=6$. Thus, $n=7$ and $p=2$ satisfy the inequality assumptions in Theorem \ref{thm:ftt2}. Since $\mathcal S_7=\mathcal K_7=\mathbf H_1$, we have ${\rm cpadd}_1^k(\mathcal S_n,\mathcal F) = {\rm cpadd}_1^k(\mathbf H_1,\mathcal F) = \mathbf H_{1+k}$ for all $k\ge0$. In particular, for every $k\in [p]$, we have ${\rm cpadd}_1^k(\mathcal S_n,\mathcal F) = \mathbf H_{1+k} = \mathcal K_{n+k} =\mathcal S_{n+k}$. Therefore, Theorem \ref{thm:ftt2} implies that $\mathbf H_k = {\rm cpadd}_1^{k-1}(\mathcal S_n,\mathcal F) =\mathcal S_{n+k-1} = \mathcal K_{6+k}$, for all $k\ge 4$. This completes the proof of (5.2.1).

\noindent (5.2.2) {\it Let $G$ be an irreducible prime graph that contains an induced $\ov{M_3^*}$. Then $G$ is $\{\mathbb L_1, \mathbb L_3\}$-free if and only if $G\in \cal N$.} 

Let $J=\ov G$ and $\mathcal F =\{\ov{\mathbb L}_1,\ov{\mathbb L}_3\}$. We prove the equivalent assertion that $J$ is $\cal F$-free if and only if $J\in\cal M$. Our proof follows the same approach as that of (5.2.1). By Theorem \ref{thm:chain-irr}, $J$ is $\cal F$-free if and only if $J\in \bigcup_{j\ge0}\mathbf J_j$, where $\mathbf J_j= {\rm cpadd}_2^j(\{M_3^*\},\cal F)$. Since $\mathbf J_0=\{M_3^*\}=\mathcal M_7$ and $\mathbf J_j = {\tt coAdd}(\mathbf J_{j-1}, \mathcal F,2,1)$ for all $j\ge1$, a computer-assisted computation (see Appendix) yields that, for each $j\in [10]$, $\mathbf J_j = \mathcal M_{7+j}$, and that the complement of every graph in $\mathbf J_j$ is irreducible. Since, by definition, $\mathcal M = \bigcup_{j\ge0}\mathcal M_{7+j}$, it remains to show that $\mathbf J_j =\mathcal M_{7+j}$ for all $j\ge11$.

Let $R=\{w_1,w_2,w_3\}$. Let $\cal A$ consist of the three $R$-rooted graphs $M_3[R\cup \{u_0\}]$, $M_3[R\cup \{v_0\}]$, and $M_3[R\cup \{u_0,v_0\}]$, and let $\cal B$ consist of the three $R$-rooted graphs $M_3[R\cup \{u_1\}]$, $M_3[R\cup \{v_1\}]$, and $M_3[R\cup \{u_1,v_1\}]$. Let $\cal S=\sigma(A,B)$. By the definitions of $\cal A, B$, and $\sigma$, we have $\cal M\subseteq S$. Moreover, for every $J= \sigma(J_0, J_1, ..., J_t)\in \cal S$, we have $J\in \cal M$ as long as $t\ge3$. Consequently, $\mathcal S_k=\mathcal M_k$ for all $k\ge 10$. For our choices of $\cal A, B$, and $\cal F$, we have $a=b=2$ and $f =6$. Thus $n=13$ and $p=4$ satisfy the inequality assumptions in Theorem \ref{thm:cpadd2}. Since $\mathcal S_n=\mathcal M_n=\mathbf J_6$, the above  computation result shows that, for all $k\in [p]$, ${\rm cpadd}_2^k(\mathcal S_n,\mathcal F) = {\rm cpadd}_2^k(\mathbf J_6,\mathcal F) = \mathbf J_{6+k}=\mathcal M_{n+k}=\mathcal S_{n+k}$, and the complement of every graph in $\mathbf J_{6+k}$ is irreducible. By Theorem \ref{thm:cpadd2}, ${\rm cpadd}_2^k(\mathcal S_n,\mathcal F) = \mathcal S_{n+k}$ for all $k>p$. Hence, for all $j\ge 11$, if we set $k=j-6>p$, we obtain $\mathbf J_j = \mathbf J_{6+k}=\mathcal S_{n+k}=\mathcal M_{7+j}$. This completes the proof of (5.2.2). 

We are now ready to prove Theorem \ref{thm:13free}. For the backward implication, let $G$ be an irreducible prime graph such that either $G\in\cal L\cup N$ or every component of $\ov G$ belongs to $\cal J\cup K$. If $G\in\cal L$, then $G$ is $\{\mathbb L_1,\mathbb L_3\}$-free, since both $\mathbb L_1$ and $\mathbb L_3$ are forbidden graphs for $\cal L$. If $G\in\cal N$, then the result follows from (5.2.2), since every graph in $\cal N$ contains an induced $\ov{M_3^*}$. It remains to consider the case in which every component $H$ of $\ov G$ belongs to $\cal J\cup K$. Since $H$ is $K_3$-free, $G$ is $\mathbb L_1$-free. To show that $G$ is $\mathbb L_3$-free, it suffices to show that $H$ is $\ov{\mathbb L}_3$-free, since $\ov{\mathbb L}_3$ is connected. If $H\in\cal J$, then $H$ is $\ov{\mathbb L}_3$-free because $\ov{\mathbb L}_3$ contains a 4-cycle while $H$ does not. If $H\in\cal K$, then $H$ is $\ov{\mathbb L}_3$-free by the first paragraph of the proof of (5.2.1). This completes the proof of the backward implication.

To prove the forward implication, let $G$ be irreducible, prime, and $\{\mathbb L_1,\mathbb L_3\}$-free. We consider two cases according to whether $\ov G$ is connected. First, let $\ov G$ be disconnected and $H$ be a component of $\ov G$. We show that $H\in\cal J\cup K$. Note that $H$ is $K_3$-free since any triangle $xyz$ in $H$ together with any vertex $w$ in another component of $\ov G$ would yield $G[w,x,y,z]\cong \mathbb L_1$, a contradiction. If the girth of $H$ is at least five, then $H\in\cal J$, and we are done. Thus, $H$ contains an induced 4-cycle. Since $H$ is a component of $\ov G$ and $G$ is prime, $H$ is coprime. Moreover, by  Observation~\ref{obs:basic}(5), $H$ is doubly connected, and hence (5.2.1) implies that $H\in\cal K$. Thus $H\in\cal J\cup K$, as required. 

Next, suppose $\ov G$ is connected. We prove that either $G\in\cal L\cup N$ or $\ov G\in\cal J\cup K$. We may assume $\ov G\succeq K_3$, since otherwise either $\ov G\not\succeq C_4$ and hence $\ov G\in\cal J$, or $\ov G\succeq C_4$ and hence $\ov G\in\cal K$ by (5.2.1). If $G\in\cal L$, we are done. Thus, $G\notin\cal L$, and so $G$ must contain an induced $\mathbb L_2, \mathbb L_4, \mathbb L_5, \mathbb L_6$, or $\mathbb L_7$. 

A vertex of a graph is {\it universal} if it is adjacent to all other vertices of the graph. Let $G^+$ denote the graph obtained from $G$ by adding a universal vertex. If $G$ has a universal vertex $v$, let $G^-=G-v$. Let $\mathbb L_{5a} = \ov{M_3^*}$. 
Let $\mathbb L_{6a}$ be obtained from $\mathbb L_6$ by adding a vertex adjacent to its unique cubic vertex, and let $\mathbb L_{6b}$ be obtained from $\mathbb L_6$ by adding a vertex adjacent to all its vertices of degree 4 or 5. We will use the following computation results (see Appendix): 
\begin{enumerate}[label=(\arabic*)]
\parskip -4pt
\item {\tt remedy}$(\mathbb L_2^-,\{\mathbb L_1, \mathbb L_3, \mathbb L_5\}) =\emptyset$.
\item {\tt remedy}$(\mathbb L_6^-,\{\mathbb L_1, \mathbb L_3, \mathbb L_5, \mathbb L_{6a},\mathbb L_{6b}\}) =\emptyset$.
\item {\tt remedy}$(\mathbb L_7^-,\{\mathbb L_1, \mathbb L_3, \mathbb L_4, \mathbb L_6\}) =\emptyset$.
\item {\tt remedy}$(\ov{\mathbb L}_4, \{\ov{\mathbb L}_1, \ov{\mathbb L}_3, \ov{\mathbb L}_5\}) =\emptyset$.
\item {\tt remedy}$(\ov{\mathbb L}_5,\{\ov{\mathbb L}_1, \ov{\mathbb L}_3, \ov{\mathbb L}_{5a}\}) = \emptyset$.
\item {\tt remedy}$(\ov{\mathbb L}_{6b},\{\ov{\mathbb L}_1, \ov{\mathbb L}_3, \ov{\mathbb L}_{6a}\})=\emptyset$.
\item {\tt conn(Add([$\ov{\mathbb L}_4$], [$\ov{\mathbb L}_1, \ov{\mathbb L}_3, \ov{\mathbb L}_5, \ov{\mathbb L_4}^+$],0,0))} $= \emptyset$
\item Let $\mathcal F=\{\ov{\mathbb L}_1,\ov{\mathbb L}_3, \ov{\mathbb L}_{5a}, \ov{\mathbb L_5}^+\}$. If {\tt x=Add([$\ov{\mathbb L}_5],\mathcal F$,0,0)}, {\tt y=Add(x,$\mathcal F$,0,0)}, {\tt z=Add(y,$\mathcal F$,0,0)}, then all graphs in {\tt z} are $K_3$-free. 
\item Let $\mathcal F=\{\ov{\mathbb L}_1,\ov{\mathbb L}_3, \ov{\mathbb L}_{6a}, \ov{\mathbb L_{6b}}^+\}$. If {\tt x=Add([$\ov{\mathbb L}_{6b}],\mathcal F$,0,0)}, {\tt y=Add(x,$\mathcal F$,0,0)}, {\tt z=Add(y,$\mathcal F$,0,0)}, then all graphs in {\tt z} are $K_3$-free. 
\end{enumerate}

Recall that $G\succeq \mathbb L_2, \mathbb L_4, \mathbb L_5, \mathbb L_6$, or $\mathbb L_7$. We show that $G\succeq \mathbb L_5$ or $\mathbb L_6$. If $G\succeq \mathbb L_2$, then $\mathbb L_2^-$  satisfies the requirements for $H$ in Corollary \ref{cor:remedy}. Hence, by Corollary \ref{cor:remedy} and (1), $G \succeq \mathbb L_5$. Similarly, by (3), $G\succeq \mathbb L_7$ implies $G \succeq \mathbb L_4$ or $\mathbb L_6$. Finally, if $G\succeq \mathbb L_4$, by Theorem \ref{thm:pathext}, $\ov G$ contains a connected $\ov{\mathbb L}_4$-add. By (7), $\ov G \succeq \ov{\mathbb L}_5$ or $\ov{\mathbb L_4}^+$. In the latter case, (4) gives $\ov G \succeq \ov{\mathbb L}_5$. Thus, $G$ contains $\mathbb L_5$ or $\mathbb L_6$.

Before presenting our final analysis, we make two key observations. First, suppose $H\subseteq \ov G$ is connected and $K_3$-free. Since $\ov G\succeq K_3$, we may fix a triangle $K$ of $\ov G$. Then $H\cup K$ is connected as $G$ is claw-free. Therefore, $H\cup K$ can be obtained from $H$ by at most three iterations of the add operation.  If we apply this observation to $H=\ov{\mathbb L}_5$, we conclude from (8) that if $G\succeq \mathbb L_5$, then $G\succeq \mathbb L_{5a}$ or $\ov G\succeq \ov{\mathbb L_5}^+$. Similarly, we conclude from (9) that if $G\succeq \mathbb L_{6b}$, then $G\succeq \mathbb L_{6a}$ or $\ov G\succeq \ov{\mathbb L_{6b}}^+$. 

Our second observation concerns $\mathbb L_{6a}$, which is a simplicial sum of $\mathbb L_6$ and $K_2$. If $G\succeq \mathbb L_{6a}$, then by Theorem \ref{thm:sumblock}, $G$ contains a graph detailed in items (1)–(7) of that theorem. By exhaustively generating all such graphs (see Appendix) we have confirmed that $G\succeq\mathbb L_{5a}$.

We now complete the proof. It suffices to show that $G\succeq \mathbb L_{5a}$, since (5.2.2) then implies $G\in\cal N$, which proves the theorem. If $G\succeq \mathbb L_5$, our first observation gives $G\succeq \mathbb L_{5a}$ or $\ov G\succeq \ov{\mathbb L_5}^+$. In the latter case, (5) gives $G\succeq \mathbb L_{5a}$, as desired. Suppose instead that $G\succeq \mathbb L_6$. By (2), $G \succeq \mathbb L_{6a}$ or $\mathbb L_{6b}$. In fact, our first observation (on $\mathbb L_{6b}$) implies that either $G\succeq \mathbb L_{6a}$ or $\ov G\succeq \ov{\mathbb L_{6b}}^+$. We further obtain from (6) that $G$ must contain $\mathbb L_{6a}$. Now our second observation implies that $G\succeq\mathbb L_{5a}$, which completes the proof of the theorem.
\hfill $\Box$\medskip

\section*{Declaration of AI Use}

 All mathematical ideas and proofs were developed independently by the authors without AI assistance. The authors used OpenAI's ChatGPT for language proofreading. The authors reviewed all suggested edits and take full responsibility for the final manuscript.

\label{refpage}

\newpage
%\section*{Appendix}
\section*{Appendix}

\subsection*{A. Functions}
\noindent\texttt{In [1]:}
\begin{lstlisting}
#all involved graphs
L1 = Graph({1:[2,3],2:[3],4:[]}).complement()
L2 = Graph({1:[2,5],2:[3],3:[4],4:[5], 6:[]}).complement()
L3 = Graph({1:[2],2:[3,5],3:[4],4:[5], 5:[6]}).complement()
L4 = Graph({1:[2],2:[3],3:[4], 5:[6]}).complement()
L5 = Graph({1:[2],2:[3],3:[4], 4:[5], 5:[6]}).complement()
L6 = Graph({1:[2],2:[3,5],3:[4],5:[6], 7:[]}).complement()
L7 = Graph({1:[2],3:[4],5:[6],7:[]}).complement()
(L2m := L2.copy()).delete_vertex(6)    ##L2-
(L6m := L6.copy()).delete_vertex(7)    ##L6-
(L7m := L7.copy()).delete_vertex(7)    ##L7-
(L6a := L6.copy()).add_edge([2,8])
(L6b := L6.copy()).add_edges([[1,8],[3,8],[4,8],[5,8],[6,8]])
m3s = Graph({0:[3,4],1:[3,4,5,6],2:[5,6],4:[5]})
L5a = m3s.complement()
(L4i := L4.copy()).add_vertex(7)      ##ading an isolated vertex to L4
(L5i := L5.copy()).add_vertex(7)      ##ading an isolated vertex to L5
(L6bi := L6b.copy()).add_vertex(9)    ##ading an isolated vertex to L6b
K332 = Graph({1:[2,6],3:[2,4],5:[2,4,6]})
\end{lstlisting}
\noindent\texttt{In [2]:}
\begin{lstlisting}
#return true if |x| = |y| and every graph in x is isomorphic to a graph in y
#note that: if graphs in x are non-isomorphic, then TRUE means x=y up to isomorphism
def iso(x,y):
    return (len(x) == len(y) and
            all(any(g.is_isomorphic(h) for h in y) for g in x))

##Returns a list containing the complement of each graph in the list G.
def co(G):
    return [g.complement() for g in G]

def conn(g):
    return [x for x in g if x.is_connected()]
    
## returns True iff G is Graph_List-free
def FreeQ(G, Graphs_List):
    return not any(G.subgraph_search(h, induced=True) for h in Graphs_List)

##returns all prime members of Graph_List
def PrimeSet(Graph_List):
    return [G for G in Graph_List if all(set(G.neighbor_iterator(u, closed=True)) != 
            set(G.neighbor_iterator(v, closed=True)) for u, v in G.edge_iterator(labels=False))]
\end{lstlisting}
\noindent\texttt{In [3]:}
\begin{lstlisting}
## returns, up to isomorphism, all Forbidden_List-free G-adds (of different versions) for G in Graph_List
## version takes three valuses 0, 1, and 2 corresponding to Add, Add1, Add2 
def Add(Graph_List, Forbidden_List, version=0, prime=True):
    new = []
    for G in Graph_List:
        V = G.vertices()
        cond = sum(G.subgraph(G.neighbors(v)).is_clique() for v in V) > 1
        s, t = (1, None) if version == 0 else \
               (1, -1) if version == 1 else \
               ((1, None) if cond else (1 + G.order(), None))
        for neighbors in Subsets(V)[s:t]:
            H = G.copy()
            u = H.add_vertex()
            H.add_edges((u, v) for v in neighbors)
            if FreeQ(H, Forbidden_List) and not any(H.is_isomorphic(X) for X in new):
                new.append(H)
    return PrimeSet(new) if prime == True else new
\end{lstlisting}
\noindent\texttt{In [4]:}
\begin{lstlisting}
## Returns up to isomorphism, all minimal F-free h-remedies
def Remedies(h, F):    
    G = []
    Vh = h.vertices()
    
    h0 = h.copy()
    u = h0.add_vertex()
    h0.add_edges((u, x) for x in Vh)

    for X in Subsets(Vh)[1:-1]:
        g = h0.copy()
        v = g.add_vertex()
        g.add_edges((v, x) for x in X)
        if FreeQ(g, F + G):
            G.append(g)

    for X in Subsets(Vh)[1:-1]:
            g = h0.copy()
            v = g.add_vertex()
            w = g.add_vertex()
            g.add_edges((v, x) for x in list(X) + [u, w])
            if FreeQ(g, F + G):
                G.append(g)
                
    return G
\end{lstlisting}

\subsection*{B. Computations for $\mathcal{M}$}
\noindent\texttt{In [1]:}
\begin{lstlisting}
#computations in the proof of (5.2.2)
load('verification_functions.sage')

#constructing all graphs in calM of order n (each consists of a base and a top)
def mn(n):
    x1 = base(n-1)
    x2 = base(n-1)
    x3 = base(n-2)
    for g in x1: g.add_edges([(0,n-1),(1,n-1),(2,n-1)]) 
    for g in x2: g.add_edges([(1,n-1)]) 
    for g in x3: g.add_edges([(0,n-2),(1,n-2),(2,n-2),(1,n-1),(n-2,n-1)]) 
    return base(n) + x1 + x2 + x3

def base(n):
    if n < 7: return []
    if n == 7: return [m3s.copy()]
    g = Graph([(0,2*i-1) for i in range(2,n//2+1)] +
         [(1,i) for i in range(3,n)] +
         [(2,2*i) for i in range(2,(n-1)//2+1)] + 
         [(2*i-1,2*i) for i in range(2,(n-1)//2+1)])
    if n % 2 == 0: return [g]
    h = g.copy(); h.delete_edge(3,4)
    return [g,h]

#check if a graph is irreducible: for a claw-free prime graph g, reducible <-> 
#(not 2-connected) or (2-connected and has a maximal clique Q of size >= 4 such that 
#g-Q has at >= 2 components and every v in Q has neighbors in <= 1 component of g-Q)
def irrQ(g):
    if g.vertex_connectivity() < 2:
        return False
    for Q in g.cliques_maximal():
        h = g.copy()
        h.delete_vertices(Q)
        C = h.connected_components()
        if len(Q) >= 4 and len(C) > 1:
            if all(sum(set(g.neighbors(v)) & set(c) 
                       != set() for c in C) <= 1 for v in Q):
                return False
    return True    
\end{lstlisting}
\noindent\texttt{In [2]:}
\begin{lstlisting}
#Verify that J_k-7 = calM_k for all k in [8,17], and all graphs are co-irreducible
F = [L1.complement(),L3.complement()]
Jk = [m3s]

for k in range(8, 18):
    Jk = co(Add(co(Jk),co(F),2,1))
    if not iso(Jk, mn(k)) or not all(irrQ(x.complement()) for x in Jk):
        print(f"A condition is violated when k = {k}")
        break
else: print(f"Conditions are verified for all k between 8 and {k}")
\end{lstlisting}

    \begin{lstlisting}[language={},frame=none,backgroundcolor=\color{white}]
Conditions are verified for all k between 8 and 17
    \end{lstlisting}

\subsection*{C. Verification}
\noindent\texttt{In [1]:}
\begin{lstlisting}
load('verification_functions.sage')
from itertools import product
\end{lstlisting}
\noindent\texttt{In [2]:}
\begin{lstlisting}
#Computations in the proof of (5.2.1)

#generate all graphs in calK with n vertices
def generateK(n):
    t = n // 2
    def K(a, b, k):
        G = graphs.CompleteBipartiteGraph(a, b)
        G.delete_edges((i, a + i) for i in range(k))
        return G
    if n % 2: return [K(t, t + 1, t)]
    return K(t, t, t), K(t, t, t - 1)

F = [graphs.CompleteGraph(3),L3.complement()]
hk = [K332]
for k in range(7, 10):
    hk = co(Add(co(hk),co(F),1,1))
    print(iso(hk, generateK(k)))
\end{lstlisting}

    \begin{lstlisting}[language={},frame=none,backgroundcolor=\color{white}]
True
True
True
    \end{lstlisting}
\noindent\texttt{In [3]:}
\begin{lstlisting}
# Computations (1)-(7) in the proof of Theorem 5.2.

sum([Remedies(L2m,[L1,L3,L5]), 
     Remedies(L6m,[L1,L3,L5,L6a,L6b]),
     Remedies(L7m,[L1,L3,L4,L6]),
     Remedies(L4.complement(),co([L1,L3,L5])), 
     Remedies(L5.complement(),co([L1,L3,L5a])),
     Remedies(L6b.complement(),co([L1,L3,L6a])),
     conn(Add(co([L4]), co([L1,L3,L5,L4i]),0,0))],[])
\end{lstlisting}
\noindent\texttt{Out [3]:}
\begin{lstlisting}[language={},frame=none,backgroundcolor=\color{white}]
[]
\end{lstlisting}
\noindent\texttt{In [4]:}
\begin{lstlisting}
# Computations (8)-(9) in the proof of Theorem 5.2.

f = co([L1, L3, L5a, L5i])
x = Add(co([L5]),f,0,0)
y = Add(x,f,0,0)
z = Add(y,f,0,0)
print([g for g in z if not FreeQ(g,[graphs.CompleteGraph(3)])])
f = co([L1, L3, L6a, L6bi])
x = Add(co([L6b]),f,0,0)
y = Add(x,f,0,0)
z = Add(y,f,0,0)
print([g for g in z if not FreeQ(g,[graphs.CompleteGraph(3)])])
\end{lstlisting}

    \begin{lstlisting}[language={},frame=none,backgroundcolor=\color{white}]
[]
[]
    \end{lstlisting}
\noindent\texttt{In [5]:}
\begin{lstlisting}
#for H = L6a, generate all graphs satisfying (1-7) of theorem 2.4 
#and confirm that no such graph is [L1, L3, L5a]-free
H, F = L6a, [L1, L3, L5a]
subsets = tuple(map(set, Subsets(H.vertices())))
s, N1, N2, R1, R2 = 2, {4, 6, 7}, {8}, {1, 3, 5}, set()
A1, A2 = R1 | N1, R2 | N2
results = []

def ext(X, Y=None, adj=False):
    G = H.copy()
    x = G.add_vertex()
    G.add_edges((x, u) for u in X)
    if Y is not None:
        y = G.add_vertex()
        G.add_edges((y, u) for u in Y)
        if adj: G.add_edge(x, y)
    return G

one_vertex_cases = [
    lambda X: (X & A1) and (N1 - X) and (s in X),
    lambda X: (X & R1) and (X & A2) and (s in X),
    lambda X: (X & A1) and (X & A2) and (s not in X)]
two_vertex_cases = [
    (False, lambda X, Y: (X & A1) and (s in X) and (Y & A1) and (s in Y)),
    (True,  lambda X, Y: (X & R1) and (s in X) and not (Y & A1)),
    (True,  lambda X, Y: (X & A1) and (s not in X) and not (Y & A1)),
    (True,  lambda X, Y: (X & A1) and not (Y & R1) and (s in Y - X))]

for condition in one_vertex_cases:
    results.append(not any(condition(X) and FreeQ(ext(X), F) for X in subsets))

for adj, condition in two_vertex_cases:
    results.append(not any(
            condition(X, Y) and FreeQ(ext(X, Y, adj), F)
            for X, Y in product(subsets, repeat=2)
        ))

print("PASS" if len(results) == 7 and all(results) else "FAIL")
\end{lstlisting}

    \begin{lstlisting}[language={},frame=none,backgroundcolor=\color{white}]
PASS
    \end{lstlisting}

\end{document}